\documentclass[11pt,reqno,natbib]{amsproc} 
	\allowdisplaybreaks
	\usepackage{graphicx}
	\usepackage{fullpage}
	\usepackage{mathptmx}
	\usepackage{enumerate}
	\usepackage{array}
	\usepackage{float}
	\floatstyle{plaintop}
	\restylefloat{table}
	\usepackage{verbatim}
	\numberwithin{equation}{section}
	\usepackage{natbib} 
	\usepackage{mathtools}
	\graphicspath{{figures/}}
	\usepackage{subfig}
	\usepackage{color}
	\usepackage{tcolorbox}
	\usepackage{booktabs}
	\usepackage{caption,setspace}
	\usepackage[debug=false, colorlinks=true, pdfstartview=FitV, 
	linkcolor=blue, citecolor=blue, urlcolor=blue]{hyperref}

	\renewcommand*{\vec}[1]{\boldsymbol{#1}}

\usepackage{morefloats}

\newlength{\drop}
\definecolor{amethyst}{rgb}{0.6, 0.4, 0.8}
\definecolor{burgundy}{rgb}{0.5, 0.0, 0.13}

\begin{document}
	
	\title{Application of the Fast Multipole Fully Coupled Poroelastic Displacement 
	Discontinuity Method to Hydraulic Fracturing Problems}
	
	\author{A.~Rezaei \and
	F.~Siddiqui \and
	G.~Bornia \and
	M.~Soliman \\
	{\small University of Houston, Houston, Texas, USA 77204--4003.\\
	 Texas Tech University, Lubbock, Texas, USA, 79409}
	}

    \begin{titlepage}
      \drop=0.1\textheight
     \centering
     \vspace*{\baselineskip}
     \rule{\textwidth}{1.6pt}\vspace*{-\baselineskip}\vspace*{2pt}
     \rule{\textwidth}{0.2pt}\\[\baselineskip]
     {\LARGE \textbf{\color{burgundy}
     	Application of the Fast Multipole Fully Coupled Poroelastic
     	Displacement Discontinuity Method to Hydraulic Fracturing
     	Problems}}\\[0.2\baselineskip]
     \rule{\textwidth}{0.4pt}\vspace*{-\baselineskip}\vspace{3.2pt}
     \rule{\textwidth}{1.6pt}\\[0.5\baselineskip]
     \scshape
     An e-print of the paper will be made 
     available on arXiv. \par
     \vspace*{0.3\baselineskip}
     Authored by \\[0.5\baselineskip]
     {\Large Ali Rezaei}, {\itshape Postdoctoral Research Associate\par}
     {\itshape Department of Petroleum Engineering \\
     	University of Houston, Houston, Texas 77204--4003 \\
     	\textbf{phone:} +1-651-245-7918, 
     	\textbf{e-mail:} arezaei2@uh.edu\par}
     \vspace*{0.75\baselineskip}
     {\Large Fahd Siddiqui}, {\itshape Postdoctoral Research Associate\par}
     {\itshape Department of Petroleum Engineering \\
     	University of Houston, Houston, Texas 77204--4003 \\
     	\textbf{phone:} +1-713-743-6103, \textbf{e-mail:} fsiddiq6@central.uh.edu 
     	\par}
     \vspace*{0.75\baselineskip}
     {\Large Giorgio Bornia}, {\itshape Associate Professor\par}
     {\itshape Department of Mathematics and Statistics \\
     	Texas Tech University, Lubbock, Texas 79409. \\ 
     	\textbf{phone:} +1-806-834-8754, \textbf{e-mail:} giorgio.bornia@ttu.edu \\
     	\textbf{website:} \url{http://www.math.ttu.edu/~gbornia/}\par}
     \vspace*{0.75\baselineskip}
     {\Large Mohammed Soliman},  {\itshape Professor\par} 
     {\itshape Department Chair and William C. Miller Endowed Chair 
     	Professor \\
     	Department of Petroleum Engineering \\
     	University of Houston, Houston, Texas 77204--4003 \\
     	\textbf{phone:} +1-713-743-6103, \textbf{e-mail:} 
     	msoliman@central.uh.edu\par}
     \vfill
     {\scshape 2019} \\
     {\small University of Houston} \par
    \end{titlepage}

	\begin{abstract}
        In this study, a  fast multipole method (FMM) is used to decrease the 
        computational time 
        of a fully-coupled poroelastic hydraulic fracture model with a controllable 
        effect on its accuracy. 
        The hydraulic fracture model is based on the poroelastic formulation of the 
        displacement discontinuity method (DDM) 
        which is a special formulation of boundary element method (BEM). 
        DDM is a powerful and efficient method for problems involving fractures.  
        However, this method becomes slow as the number of temporal, or spatial 
        elements increases,
        or necessary details such as poroelasticity, that makes the solution 
        history-dependent, are added to the model. 
        FMM is a technique to expedite matrix-vector multiplications
        within a controllable error without forming the matrix explicitly.
        Fully-coupled poroelastic formulation of DDM involves the multiplication of 
        a dense matrix with a vector in several places.
        A crucial modification to DDM is suggested in two places in the algorithm 
        to leverage the speed efficiency of FMM for carrying out these 
        multiplications.
        The first modification is in the time-marching scheme,
        which accounts for the solution of previous time steps to compute the 
        current time step.
        The second modification is in the generalized minimal residual method 
        (GMRES) 
        to iteratively solve for the problem unknowns.
        
        Several examples are provided to show the efficiency of the proposed 
        approach 
        in problems with large degrees of freedom (in time and space).
        Examples include hydraulic fracturing of a horizontal well
        and randomly distributed pressurized fractures at different orientations 
        with respect to horizontal stresses. 
        The results are compared to the conventional DDM in terms of computational 
        processing time and accuracy.
        It is demonstrated that for the case of 20000 constant spatial elements and 
        a single temporal element, 
        FMM may decrease the computation time by up to 70 times 
        with a relative error less than $ 4 $\% for the hydraulic fracture example,
        and less than $ 0.5 $\% for the case of randomly distributed pressurized 
        fractures. 
        The solution of tip displacements using both methods are then used to 
        compare the computation of stress intensity factors (SIF) in mode I and II, 
        which are needed for fracture propagation. 
        The error of SIF calculation using the proposed modification was also found 
        to be negligible. 
        Therefore, this method will not affect the estimation of the fracture 
        propagation direction.
        Accordingly, the proposed algorithm may be used for fracture propagation 
        studies 
        while substantially reducing the processing time.
		\keywords{}
	\end{abstract}

    \footnotetext{\emph{Key words} Fast Multipole Method; Boundary Element Method; 
    Displacement
Discontinuity Method; Hydraulic fracturing; Poroelasticity}
	
	\maketitle

	\section{Introduction}
	The boundary element method (BEM) \citep{jaswon1963integral, rizzo1967integral, 
	banerjee1981boundary, aliabadi1991boundary, cruse2012boundary}
	is a well-known and efficient method for solving problems 
	with high volume-to-surface ratio.
	An example of such a problem is the hydraulic fracturing process, 
	which is used in the oil and gas industry to increase the production of 
	hydrocarbon from tight reservoirs. 
	For a comprehensive review on advances of BEM one may refer to 
	\cite{liu2011recent}.
	A special formulation of the boundary element method, known as displacement 
	discontinuity method \citep{crouch1976solution}, 
	is extensively used to study hydraulic fracturing problems 
	\citep{curran1987displacement, carvalho1991poroelastic, 
	detournay1989poroelasticity, 
		ghassemi2013three, peirce2014robustness, safari20143d, wu2015simultaneous, 
		Rezaei_Journal2018}. 
	The host medium of hydraulic fractures is poroelastic and contains 
	discontinuities. 
	Therefore, reliably modeling the behavior of hydraulic fractures 
	and their complex interaction with the surrounding environment 
	requires coupling between various phenomena, such as flow of fracturing fluid 
	inside the fracture, 
	flow of the fluids in the poroelastic rock, 
	deformation of the porous rock, 
	fluid leak off from fracture, 
	and hydraulic fracture propagation. 
	However, including the necessary details 
	such as strong coupling between pore pressure and rock displacement in the BEM 
	formulation 
	makes the method computationally inefficient. 
	This is because of the requirement for discretization both space and time
	that causes the solution to become history-dependent. A summary on other 
	challenges that are involved in hydraulic fracturing problems are presented by 
	\cite{peirce2016implicit}. 
	
	One main difference between hydraulic fracture models is the approach in 
	handling fracturing fluid leak-off. The calculations of leak-off, width, and 
	length using hydraulic fracture models 
	may be categorized into three groups based on complexity of the 
	relationship between fracturing fluid diffusion and rock deformation 
	\citep{vandamme1990poroelasticity}. 
	The first approach is called \textit{uncoupled model}. 
	In this category, the main assumption is that the rock is linearly elastic. 
	Therefore, no fluid flow inside the rock pore space is assumed 
	and leak-off is calculated by a one dimensional Carter's model 
	\citep{howard1970hydraulic}.
	The second category consists of partially coupled models.
	In these models, stresses and displacements are still based on the theory of 
	elasticity. 
	In these models, the effect of leak-off is considered by the linear diffusion 
	law.
	Also, the concept of back stress \citep{cleary1980comprehensive} is used
	to consider the effect of pore pressure in these models. 
	The third category, which is used in this study, belongs to the class of 
	fully-coupled models. 
	These models are based on \cite{biot1941general} theory of poroelasticity. 
	In these models, a full range 
	of coupled diffusion-deformation are considered.
	
	The fast multipole method (FMM) 
	\citep{rokhlin1985rapid,greengard1987fast,Greengard:1987:REP:913529} 
	is one of the top ten scientific computing algorithms that were developed in 
	twentieth century \citep{liu2006fast}. 
	FMM is a technique to expedite matrix-vector multiplications within a 
	controllable error without forming the matrix explicitly. 
	This method was initially introduced to solve a reciprocal function of the 
	distance between two points 
	using Legendre polynomials \citep{greengard1987fast}. 
	Since then, different fast summation techniques have been developed 
	\citep[e.g.][]{alpert1993wavelet, gimbutas2001coulomb, gimbutas2003generalized, 
	ying2004kernel,cheng2005compression, dahmen2006compression, 
	martinsson2007accelerated, wengmodeling2015}. 
	A review of the progress of the fast multiple method is presented by 
	\cite{nishimura2002fast}, \cite{liu2009fast}, and \cite{yokota2016fast}. 
	The main reason for developing different FMM techniques is that in the original 
	method,
	an analytical expansion of the kernel was required limiting its applications to 
	more specific cases.
	In order to overcome this problem, a set of fast multipole techniques 
	were developed that rely only on the numerical values of the kernel function. 
	\cite{fong2009black} introduced a kernel-independent fast multipole method. 
	Their approach is useful for the kernels for which analytical expansions are 
	not known. 
	It uses Chebyshev polynomials for expanding the kernel functions.

	To overcome the deficiency of the fully-coupled version of DDM in large 
	problems different approaches may be taken. 
	\cite{cheng2016rapid} suggested a method to overcome this deficiency by 
	approximating the solution using an analytical formulation. 
	In another approach, FMM may be used to expedite matrix-vector multiplications. 
	Several studies have applied different fast multipole techniques to BEM 
	in problems involving fractures 
	\citep[e.g.][]{ nishimura1999fast, helsing2000fast, lai2003fast, otani2008fmm, 
	wang2011fast, guo2014fast, liu2017modeling}.
	Most of these problems either utilized an analytical kernel expansion 
	\citep[][]{yoshida2001application,wang2005fast, liu2009fast,liu2017modeling}, 
	or applied the numerical kernel expansion to elastic problems 
	\citep[][]{farmahini2016simulation, verde2015modeling}, 
	partially coupled problems \citep[][]{verde2016large}, or fully coupled 
	poroelastic media \citep[][]{schanz2018fast}. 
	\cite{morris2000efficient} utilized a modified FMM with the DDM formulation to 
	study a rock sample failure Brizilian test.
	\cite{peirce1995spectral} reduced the BEM cost to $ O(N^2\log N) $ operations 
	using the Spectral Multipole Method (SMM). \cite{liu2017modeling} studied the 
	propagation of multiple fractures in an elastic solid using a dual boundary 
	integral equation (BIE) and FMM. The black-box fast multipole method (bbFMM) 
	was applied to the elastic formulation of the displacement discontinuity method 
	\citep{verde2013efficient,verde2013fracture,verde2015fast,verde2015modeling,farmahinianalysis,farmahini2016simulation}
	 to study problems such as simulation of micro-seismicity in response to 
	injection/extraction in fracture networks.  
	\cite{verde2016large} applied the bbFMM to a partially-coupled formulation of 
	DDM to study permeability variation in fractured reservoirs.

	In this study, the FMM described by \cite{fong2009black} is implemented into a 
	conventional fully-poroelastic DDM to solve the problems of hydraulic 
	fractures. The fully-coupled formulation of DDM that is used in this study and 
	application of the developed model is the distinction from studies presented by 
	\cite{verde2015fast, verde2015modeling, 
	verde2016large,farmahini2016simulation}. Several examples are provided to show 
	the efficiency of the proposed approach 
	in problems with large degrees of freedom (in time and space). 
	Examples include hydraulic fracturing of a horizontal well
	and randomly distributed pressurized fractures at different orientations 
	with respect to horizontal stresses. 
	Results are compared to the conventional DDM in terms of computational 
	processing time and accuracy. 
	Details of the procedure and implementation of the method 
	to hydraulic fracture problems will be presented in the following sections.
	The outline of the paper is as follows. 
	The poroelastic displacement discontinuity method is reviewed in Section 
	\ref{sec:pddm}. 
	Then, a discussion on the calculation of stress intensity factors using DDM is 
	presented in Section \ref{sec:sif}. 
	Next, the FMM and solution procedure using the proposed approach are explained 
	in Section \ref{sec:fmm}. 
	Finally, two examples of applications of the proposed approach are explained in 
	Section \ref{sec:results}.

	\section{Poroelastic Displacement Discontinuity Method (PDDM)} \label{sec:pddm}
	
\begin{figure}[H]
		\centering
		\includegraphics[scale=0.6]{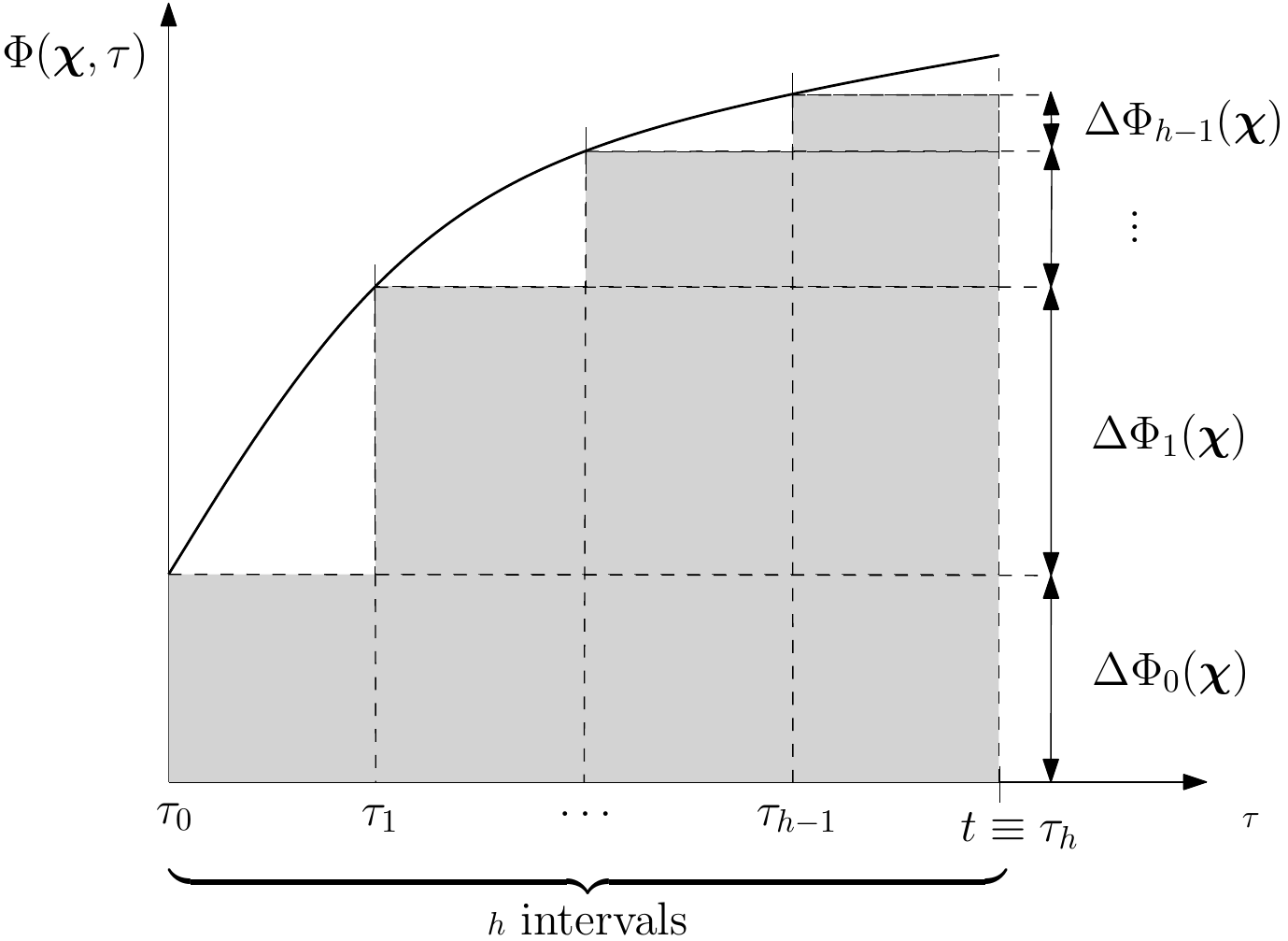}
		\caption{Time marching scheme} 
		\label{fig:time_marching}	
	\end{figure}

	The displacement discontinuity method 
	(DDM) \citep{crouch1976solution} is a special formulation of BEM 
	in continuum media involving a fracture. It is an indirect boundary element 
	method (BEM), 
	formulated for media containing cracks, where a discontinuity exists in 
	displacements. 
	The initial formulation of DDM was based 
	on purely elastic medium and may be derived from dislocation theory 
	\citep{bobet2005stress}. 
	\cite{liu2014revisit} explicitly showed that DDM and BEM are equivalent for 
	fracture problems.  
	Green's functions (fundamental solutions) are key elements of any BEM 
	formulation.
	Fundamental solutions of poroelastic medium were derived by 
	\cite{cleary1977fundamental} 
	based on the governing equations of the theory of poroelasticity 
	\citep{biot1941general}. 
	\cite{curran1987displacement} and \cite{detournay1987poroelastic} 
	used the fundamental solutions for poroelastic media 
	to develop the solutions of poroelastic media using DDM. 
	The poroelastic displacement discontinuity method 
	allows for the calculation of changes in pore pressure, stress, and 
	displacement over time.
	The boundary integral equations (BIE) relating stress and pore pressure 
	along the fracture and in the rock medium 
	to displacements and fluid leak-off may be written as
	\begin{align}
	&
	\begin{aligned}  \label{eq:int_sigs}
	\sigma_{ij}(\vec{x},t) = 
	& 
	\int_{0}^{t} \int_{\Gamma}
	Q_{ik}(\vec{\chi}) Q_{jl}(\vec{\chi})  S_{s,kl}(\vec{x},\vec{\chi};t-\tau) 
	D_{s}(\vec{\chi},\tau) 
	d\Gamma(\vec{\chi}) d\tau 
	\\
	+ & \int_{0}^{t} \int_{\Gamma}
	Q_{ik}(\vec{\chi}) Q_{jl}(\vec{\chi})  S_{n,kl}(\vec{x},\vec{\chi};t-\tau) 
	D_{n}(\vec{\chi},\tau)
	d\Gamma(\vec{\chi}) d\tau 
	\\
	+ & \int_{0}^{t} \int_{\Gamma}
	Q_{ik}(\vec{\chi}) Q_{jl}(\vec{\chi}) S_{q,kl}(\vec{x},\vec{\chi};t-\tau) 
	D_{q}(\vec{\chi},\tau)  
	d\Gamma(\vec{\chi}) d\tau \,,
	\end{aligned} \\
	&
	\begin{aligned}
	p(\vec{x},t) = 
	& 
	\int_{0}^{t}\int_{\Gamma} 
	P_{s}(\vec{x},\vec{\chi};t-\tau) D_{s}(\vec{\chi},\tau)
	d\Gamma(\vec{\chi}) d\tau 
	\\
	+ &
	\int_{0}^{t}\int_{\Gamma} 
	P_{n}(\vec{x},\vec{\chi};t-\tau) D_{n}(\vec{\chi},\tau)
	d\Gamma(\vec{\chi}) d\tau 
	\\
	+ &
	\int_{0}^{t}\int_{\Gamma} 
	P_{q}(\vec{x},\vec{\chi};t-\tau) D_{q}(\vec{\chi},\tau)
	d\Gamma(\vec{\chi}) d\tau 
	\,.
	\end{aligned}
	\label{eq:int_pp}
	\end{align}
	In Equations \eqref{eq:int_sigs} - \eqref{eq:int_pp}, $\sigma_n$, $\sigma_s$, 
	$p$ 
	are 
	normal stress, shear stress, and pore pressure respectively. 
	Also, $ D_s $, $ D_n $, and $ D_q $ are shear, normal, and fluid loss. 
	The integral equations \eqref{eq:int_sigs} - \eqref{eq:int_pp} and their 
	fundamental solutions 
	assume instantaneous impulses. 
	Therefore, an integration over time is required to get the fundamental 
	solutions of continuous impulses. 
	In order to account for the temporal part of the integral in Equations 
	\eqref{eq:int_sigs} - \eqref{eq:int_pp}, 
	different techniques may be utilized. 
	The time marching technique is used in this study to solve the problem at 
	successive time intervals.
	This technique leads to systems of linear equations 
	which are simultaneous in space but successive in time 
	\citep{banerjee1981boundary}. 
	Figure \ref{fig:time_marching} shows a schematic of the time marching process.
	
	The discretization of the Equations \eqref{eq:int_sigs} - \eqref{eq:int_pp} 
	using constant spatial and constant temporal elements
	may be written as 
	\begin{equation}
	\begin{split}
	\sum\limits_{\lambda=1}^N 
	A_{xx}^{\beta \lambda}D_{s}^{\lambda,h} +
	\sum\limits_{\lambda=1}^N 
	A_{xy}^{\beta \lambda}D_{n}^{\lambda,h} +
	\sum\limits_{\lambda=1}^N 
	A_{xq}^{\beta \lambda}D_{q}^{\lambda,h} =
	\\\sigma_{s}^{h}(x^\beta,t) - 
	\sum\limits_{\eta=0}^{h-1}\sum\limits_{\lambda=1}^N
	\bigg(
	A_{xx}^{\beta \lambda,\eta}D_{s}^{\lambda,\eta} + 
	A_{xy}^{\beta \lambda,\eta}D_{n}^{\lambda,\eta} + 
	A_{xq}^{\beta \lambda,\eta}D_{q}^{\lambda,\eta} \bigg), 
	\end{split}
	\label{eq:final_sigs}
	\end{equation}
	\begin{equation}
	\begin{split}
	\sum\limits_{\lambda=1}^N 
	A_{yx}^{\beta \lambda}D_{s}^{\lambda,h} +
	\sum\limits_{\lambda=1}^N 
	A_{yy}^{\beta \lambda}D_{n}^{\lambda,h} +
	\sum\limits_{\lambda=1}^N 
	A_{yq}^{\beta \lambda}D_{q}^{\lambda,h} =
	\\\sigma_{n}^{h}(x^\beta,t) - 
	\sum\limits_{\eta=0}^{h-1}\sum\limits_{\lambda=1}^N
	\bigg(
	A_{yx}^{\beta \lambda,\eta}D_{s}^{\lambda,\eta} + 
	A_{yy}^{\beta \lambda,\eta}D_{n}^{\lambda,\eta} + 
	A_{yq}^{\beta \lambda,\eta}D_{q}^{\lambda,\eta} \bigg), 
	\end{split}
	\end{equation}	
	\begin{equation}
	\begin{split}
	\sum\limits_{\lambda=1}^N 
	A_{px}^{\beta \lambda}D_{s}^{\lambda,h} +
	\sum\limits_{\lambda=1}^N 
	A_{py}^{\beta \lambda}D_{n}^{\lambda,h} +
	\sum\limits_{\lambda=1}^N 
	A_{pq}^{\beta \lambda}D_{q}^{\lambda,h} =
	\\p_p^{h}(x^\beta,t) - 
	\sum\limits_{\eta=0}^{h-1}\sum\limits_{\lambda=1}^N
	\bigg(
	A_{px}^{\beta \lambda,\eta}D_{s}^{\lambda,\eta} + 
	A_{py}^{\beta \lambda,\eta}D_{n}^{\lambda,\eta} + 
	A_{pq}^{\beta \lambda,\eta}D_{q}^{\lambda,\eta} \bigg). 
	\end{split}
	\label{eq:final_pp}
	\end{equation}
	In Equations \eqref{eq:final_sigs}-\eqref{eq:final_pp}, 
	$ A_{ij} $ are the coefficients relating the displacement discontinuities and 
	fluid sources 
	to shear stress, normal stress and pore pressure 
	\citep{carvalho1991poroelastic}. 
	For example, $ A_{xx} $ is the shear stress that is induced on the observation 
	point 
	from a unit shear displacement discontinuity at the influencing point. 
	In general, 
	the fracturing fluid pressure is known at the boundary of a hydraulic 
	fracturing problem (i.e. fracture surface), 
	and the fracture surface displacements and flow discontinuity 
	are the unknowns. 
	Therefore, Equations \eqref{eq:final_sigs}-\eqref{eq:final_pp} 
	form a set of $ 3 N $ linear equations 
	which may be solved for $ 3 N $ unknowns namely $\sigma_s$, $\sigma_n$, and $p$ 
	for a single time step.
	As it may be seen, totally nine matrix-vector multiplications exist on the 
	right hand side 
	of Equations \eqref{eq:final_sigs}-\eqref{eq:final_pp}. 
	Moreover, for any extra time step that is added to the problem, 
	nine matrix-vector multiplications will be added to the computation. 
	Furthermore, using an iterative solver like GMRES requires nine matrix-vector 
	multiplications at each iteration. 
	These operations increase the computational time especially with a large number 
	of spatial or temporal elements.
	Hence, to improve the computational efficiency, a fast multipole method is 
	implemented in the PDDM algorithm. \cite{Rezaei_ARMA2017} applied the bbFMM to 
	the double summation on the right hand-side of Equations \eqref{eq:final_sigs} 
	- \eqref{eq:final_pp}. In this study, bbFMM is applied on both double summation 
	and GMRES iterative solver. 
	Details of FMM and its implementation are discussed in the following sections 
	of the paper.
	
	\subsection{Calculation of Stress Intensity Factors} \label{sec:sif}
	
	Two classes of parameters are required 
	in order to determine if the fracture propagation occurs. 
	The first parameter is the fracture toughness $ K_{Ic} $,
	which measures the ability of a material containing a fracture to resist a load 
	and is determined experimentally. 
	The second set of parameters is given by the stress intensity factors (SIFs) 
	\citep{irwin1957analysis}, 
	which are a function of the fracture length and of the stress applied on the 
	surface of the fracture.
	
	Following the principles of linear elastic fracture mechanics (LEFM),
	a fracture may propagate according to three different modes: 
	opening or tensile (mode I), 
	plane shearing or sliding (mode II), 
	out-of-plane shearing or tearing (mode III).
	There are three stress intensity factors associated with these different 
	loading modes, 
	known as $ K_I $, $ K_{II} $ and $ K_{III} $, respectively. 
	Among these three, only $ K_{I} $ and $ K_{II} $ are considered in 2D cases 
	because they don't involve the third dimension. 
	The calculation of the stress intensity factors
	plays a critical role in mixed mode fracture propagation criteria.
	SIFs are usually calculated at the tip of the fracture. 
	\cite{olson1990fracture} empirically introduced a relationship for 
	calculating the SIFs using $ D_n $ and $ D_s $ of the fracture tip element as
	\begin{equation}
	K_{I}   = 0.806 \frac{E}{4(1-\nu^2)} \sqrt{\frac{\pi}{2a}} D_n \,,\quad
	K_{II}  = 0.806 \frac{E}{4(1-\nu^2)} \sqrt{\frac{\pi}{2a}} D_s  \,.
	\label{eq:SIF}
	\end{equation}
	In this paper, our aim is to compare the accuracy of the calculated SIFs 
	using either a conventional or a fast multipole version of a fully poroelastic 
	DDM model.

	\section{Fast Multipole Method (FMM)} \label{sec:fmm}
	
	The Fast Multipole Method (FMM) is a method to efficiently calculate 
	matrix-vector products of the type
	\begin{align} \label{eq:fmm}
	f(\vec{x}_i) = \sum_{j=1}^{N} G(\vec{x}_i,\vec{y}_j) \, b(\vec{y}_j) \,,\quad i 
	= 1, \ldots, N  \,,
	\end{align}
	using $ O(N) $ operations and a controllable error. 
	Here, the kernel function $ G $ 
	relates the \textit{influenced} points  $ \vec{x}_i $ 
	with the \textit{influencing} points $ \vec{y}_j $,  
	$ b(\vec{y}_j) $ are the \textit{charges} at the influencing points, 
	and $ f(\vec{x}_i) $ is the \textit{potential} at the influenced points.
	The main idea of FMM is to speed up the evaluation of the summation in 
	\eqref{eq:fmm}  
	by using a \textit{hierarchical tree decomposition} of the domain points.
	Thanks to this subdivision, a fast approximation of
	the kernel function can be introduced for large distances 
	between the influencing and the influenced points,
	while the direct multiplication in \eqref{eq:fmm} is used
	for points that are not well-separated.

	\subsection{Hierarchical Tree Decomposition}
	
	\begin{figure}[t]
		\centering
		\subfloat[Level 0]{
			\includegraphics[width=0.25\textwidth]{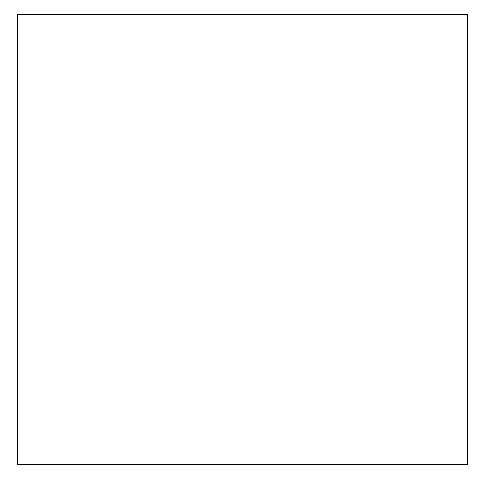}
			\label{fig:level_0}
		}
		\subfloat[Level 1]{
			\includegraphics[width=0.25\textwidth]{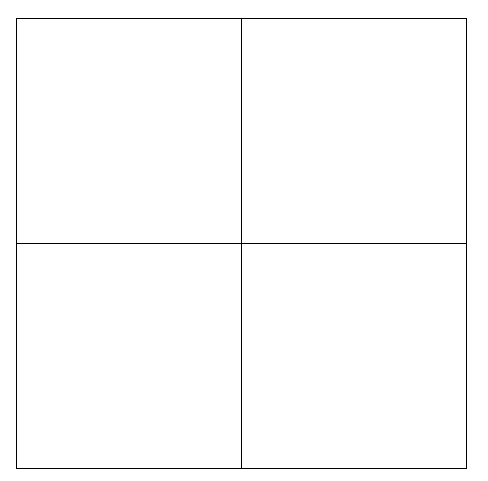}
			\label{fig:level_1}
		}	
	    \quad
		\subfloat[Level 2]{
			\includegraphics[width=0.25\textwidth]{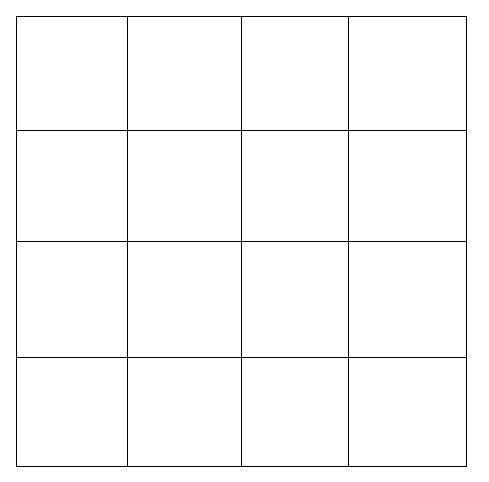}
			\label{fig:level_2}
		}	
		\subfloat{
			\includegraphics[width=0.15\textwidth]{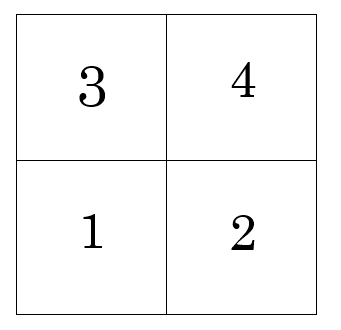}
			\label{fig:indexing}
		}	
		\caption{Hierarchical tree decomposition; (a) level 0 square, (b) level 1 
		square, (c) level 2 square, and (d) indexing convention }
		\label{fig:hei_decompose}
	\end{figure}

	The hierarchical tree decomposition 
	is the first step towards any fast multipole method. 
	Figure \ref{fig:hei_decompose} shows how this decomposition is performed in two 
	dimensions,
	in which case it is also referred to as \textit{quad-tree} structure.
	The square that covers the entire domain of the problem is called \textit{level 
	0} square (Figure \ref{fig:level_0})
	and it is first divided into 4 \textit{child} squares (Figure \ref{fig:level_1})
	that define level 1. 
	At any stage, any square that has more points than a prescribed number is 
	recursively divided into 4 child squares
	and a new level is generated.
	Otherwise, the square is defined a \textit{leaf cell}. 
	In particular, a leaf cell that does not contain any points is called a 
	\textit{zero cell}.
	The process of cell subdivision stops when reaching a leaf cell.
	When new child cells are created, new indices must be assigned. 
	Our convention is that the cell indexing increases from left to right 
	and from bottom to top, as shown in Figure \ref{fig:indexing}. 
	For example, the bottom left square and top right square at level 1 are the 
	children number 1 and 4 of the level 0 square, respectively. 
	It is worth mentioning that the indexing of cells is unrelated to the global 
	numbering of the domain points. 
	
	\begin{figure}[t]
		\centering
		\includegraphics[scale=0.35]{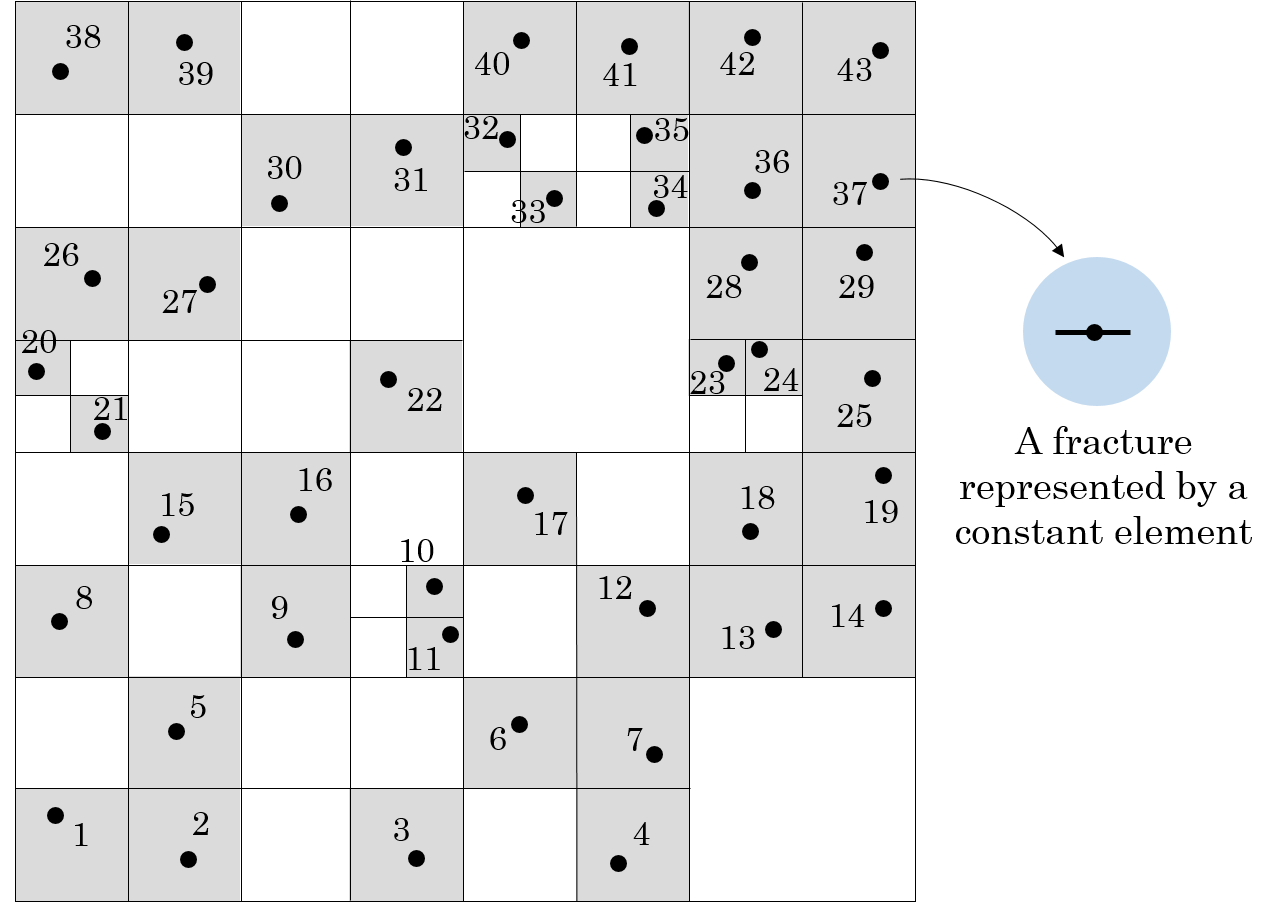}  \\ \vspace{1em}
		\includegraphics[scale=0.18]{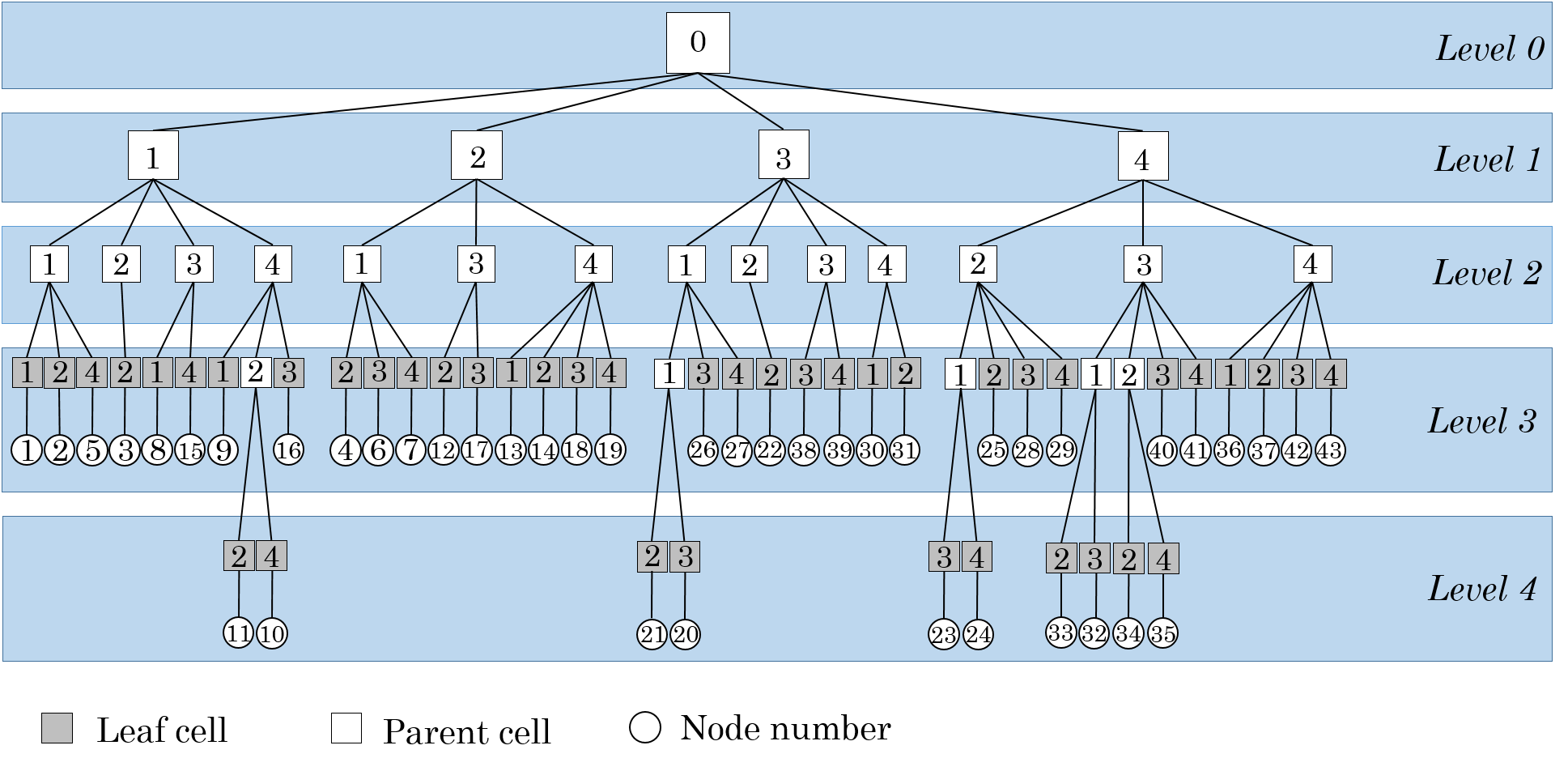}
		\caption{Top: an example with 43 randomly distributed points.
			Bottom: the associated hierarchical quad-tree structure.} 
		\label{fig:fmm_tree}	
	\end{figure}
	
	To better explain the tree structure, 
	Figure \ref{fig:fmm_tree} at the top shows an example of a domain with 43 
	randomly distributed points. 
	In the PDDM method, each of them represents the center point of a fracture 
	element.
	The hierarchical tree for this example is shown in Figure \ref{fig:fmm_tree}. 
	It summarizes the cell subdivisions
	as well as the relationships between cells and contained points.
	For example, child 2 and child 4 of level 0 have only three children each, 
	because one of their children is a zero cell.
	Each point can be traced back to level 0 using the tree. 
	For example, point number 21 belongs to the second child of the first child of 
	the first child of the third child of the level 0 cell
	(i.e., the entire domain).
	
	\begin{figure}[t]
		\centering
		\includegraphics[scale=0.5]{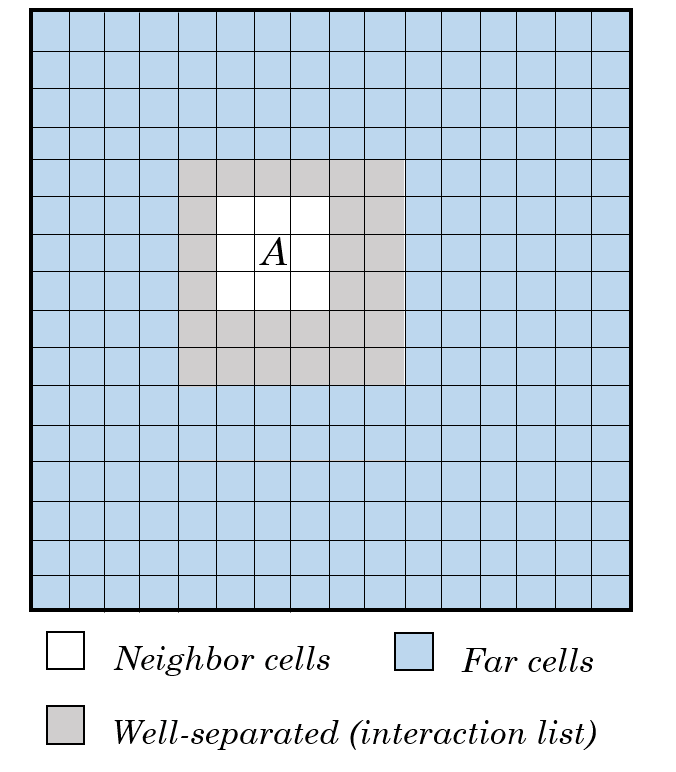}
		\caption{Three groups of cells associated with cell $A$} 
		\label{fig:fmm_grouping}	
	\end{figure} 
	For each cell, three groups may be defined by using the tree structure:
	\begin{itemize}
		\item the list of \textit{neighbor cells}, which are the cells at any level 
		that have at least one common vertex with the given cell;
		\item the \textit{interaction list} of the cell, which is made of all cells 
		that are \textit{well-separated} with respect to the given cell 
		(two cells are said to be well-separated if they are not neighbors at the 
		same level, but their parents are neighbors); 
		\item the list of \textit{far cells}, which are all the remaining cells in 
		the domain.
	\end{itemize}
	Figure \ref{fig:fmm_grouping} shows an example of neighbor, well-separated, and 
	far cells for cell $ A $.
	In the example of Figure \ref{fig:fmm_tree}, 
	the cells containing points 1 and 2 are neighbors, 
	the cells containing points 1 and 3 are in the interaction lists of each other, 
	and the cells containing point 1 and 43 are far cells.
	It should be noted that after the hierarchical tree structure is constructed, 
	it will not change unless there is a change in the problem geometry such as 
	fracture propagation. 
	After its construction, the hierarchical tree can be used to calculate the 
	potential at each point.

	\subsection{The Black-Box Fast Multipole Method}
	
	Different versions of the fast multipole method can be defined.
	We may distinguish them with respect to the way in which the kernel function is 
	approximated.
	On one hand, we may consider expanding the kernel in terms of some analytical 
	expansion
	(e.g. spherical harmonics expansion), and truncating such series.
	On the other, we may approximate the kernel through interpolation with respect 
	to some basis.
	\cite{fong2009black} introduced a version of the fast multipole method, 
	called the \textit{black-box Fast Multipole Method}, 
	which belongs to the class of interpolation-based fast multipole methods.
	This method has several advantages 
	compared to the ones based on an analytical kernel expansion.
	First of all, the method is \textit{black-box} or \textit{kernel-independent} 
	in the sense that 
	an analytical expansion of the kernel need not be known.
	Since the method is based on interpolation, only the evaluation of the kernel 
	at certain points is required.
	Moreover, the interpolation is performed using Chebyshev polynomials, which 
	have several advantages 
	such as uniform convergence and near minimax approximation 
	\citep{fong2009black}. 
	
	\subsubsection{Kernel approximation with Chebyshev interpolation}

	For the sake of simplicity, let us describe the method in a one-dimensional 
	domain.
	In this case, the corresponding two-dimensional kernel function $ G(x,y) $ in 
	Equation \eqref{eq:fmm} 
	is approximated by a \textit{low-rank approximation} $ \widetilde{G}(x,y) $ 
	using polynomials $ u_l $ and $ v_l $, for $ l = 1,2,\ldots,n $,
	\begin{align}   \label{eq:low_rank}
	\widetilde{G}(x,y) = \sum_{l=1}^{n} u_l(x) v_l(y) \,.
	\end{align}
	Combining Equations \eqref{eq:low_rank} and \eqref{eq:fmm}, 
	a fast summation scheme for an approximation $ \widetilde{f} $ of the potential 
	$ f $ is given by
	\begin{align} 	\label{eq:fast_sum}
	\widetilde{f}(x_i) = \sum_{l=1}^{n} u_l(x_i) \sum_{j=1}^{N} v_l(y_j) b(y_j).
	\end{align}
	Hence, first we transform the source charges 
	using the second summation in the right-hand side of \eqref{eq:fast_sum}. 
	Then, we calculate the influence at each observation point 
	using the first summation.
	While the computational complexity of Equation \eqref{eq:fmm} is  $ O( N^2 ) $,
	a significant reduction is given by Equation \eqref{eq:fast_sum} 
	since in this case the complexity is $ O ( 2 n N ) $,
	which is especially important when $ n \ll N $.

	Following \cite{fong2009black}, we choose Chebyshev polynomial interpolation 
	for the low-rank kernel approximation \eqref{eq:low_rank}. 
	Given a function $ g(x) $, the interpolating polynomial $ p_{n-1}(x) $ of 
	degree $ n - 1 $  
	using Chebyshev polynomials $ T_k $ can be written in the form (see 
	\cite{fong2009black})
	\begin{align}
	p_{n-1}(x) & = \sum_{l=1}^{n} g(\bar{x}_l) \, S_n(\bar{x}_l,x) \,, \text{ where 
	} \\
	S_n(x,y) & = \dfrac{1}{n} + \dfrac{2}{n} \sum_{k=1}^{n-1} T_k(x) T_k(y) \,,
	\end{align}
	and $ \bar{x}_l $ are the roots of $ T_n $.
	Using this interpolating polynomial for both $ u_l $ and $ v_l $ in the 
	low-rank approximation \eqref{eq:low_rank},
	the approximated kernel function $ \widetilde{G} $ becomes
	\begin{align}
	\widetilde{G}({x},{y}) = \sum_{l=1}^{n} \sum_{m=1}^{n} 
	G({\bar{x}_l},{\bar{y}_m}) S_n({\bar{x}_l},{x}) S_n({\bar{y}_m},{y})
	\end{align}
	Substituting this approximation into Equation \eqref{eq:fmm}, 
	the fast summation method reads
	\begin{align} 	\label{eq:fmm_2d}
	\widetilde{f}(x_i) = \sum_{j=1}^{N} \widetilde{G}(x_i,y_j) b(y_j) 
	& = \sum_{j=1}^{N} \bigg[ \sum_{l=1}^{n} \sum_{m = 1}^{n} 
	G({\bar{x}_l},{\bar{y}_m}) S_n({\bar{x}_l},x_i) S_n({\bar{y}_m},y_j) \bigg] 
	b(y_j) \\
	& = \sum_{l=1}^{n} S_n({\bar{x}_l},x_i) \sum_{m=1}^{n} 
	G({\bar{x}_l},{\bar{y}_m}) \sum_{j=1}^{N} b(y_j) S_n({\bar{y}_m},y_j). 
	\label{fastsum_cheb}
	\end{align} 
	First, one computes the weights at the Chebyshev nodes (summation over $ j $ in 
	\eqref{fastsum_cheb}), 
	then the approximated function $ \widetilde{f} $ is computed at the Chebyshev 
	nodes (summation over $ m $), 
	and finally $ \widetilde{f} $ is calculated at the influenced points by 
	interpolation (summation over $ l $). 
	
	\subsubsection{Upward Pass and Downward Pass} 
	
	
	Using the hierarchical tree decomposition, 
	clusters of particles are determined at different levels.  
	The interaction between clusters that are well-separated are calculated using 
	FMM,
	and the interaction between the clusters that are not well-separated are 
	calculated
	using direct matrix-vector multiplication. 
	Hence, two main steps have to be taken in order to compute potentials at each 
	point:
	the \textit{upward pass} and the \textit{downward pass}. 
	
	\begin{figure}[t]
		\centering
		\includegraphics[scale=0.45]{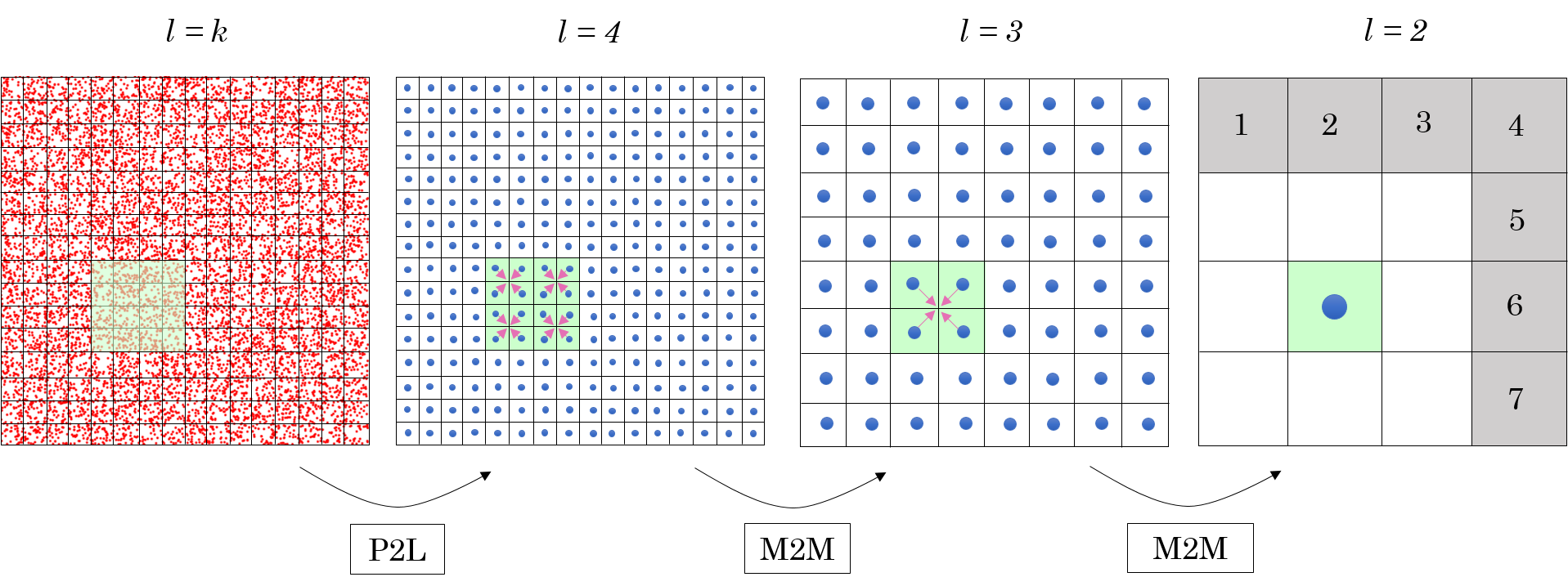}
		\caption{Upward pass}
		\label{fig:UpwardPass}
	\end{figure}
	The algorithm starts with the \textit{upward pass} (Figure 
	\ref{fig:UpwardPass}),
	whose purpose is to construct weights at all levels for a subsequent evaluation 
	of the potential.
	For each cell at all levels, 
	the weights $ W_m^I $ at the Chebyshev nodes $ \bar{y}^I_m $ are calculated as
	(here we denote the set of child cells of a given cell $ I $ as $ 
	\mathrm{Children}(I) $)
	\begin{align}
	( l = k )                        \quad  W_m^I & = \sum_{y_j \in I} b(y_j) 
	S_n(\bar{y}^I_m,y_j), \quad m = 1,...,n \,, \label{p2l} \\
	( l = k - 1, k-2, \ldots, 1, 0 ) \quad  W_m^I & = \sum_{J \in 
	\mathrm{Children}(I)}
	\sum_{m'} W_{m'}^{J} S_n(\bar{y}^I_m,\bar{y}_{m'}^{J}), \quad m = 1,...,n \,. 
	\label{m2m}
	\end{align}
	The interpolation at the finest level in \eqref{p2l} is called 
	\textit{point-to-local} (P2L) translation.
	All other levels are calculated recursively by moving through the tree upwards, 
	from finer to coarser levels.
	Equation \eqref{m2m} is referred to as \textit{multipole-to-multipole (M2M) 
	translation} and consists 
	of an interpolation of the weights from the child cells.
	
	\begin{figure}[t]
		\centering
		\includegraphics[scale=0.45]{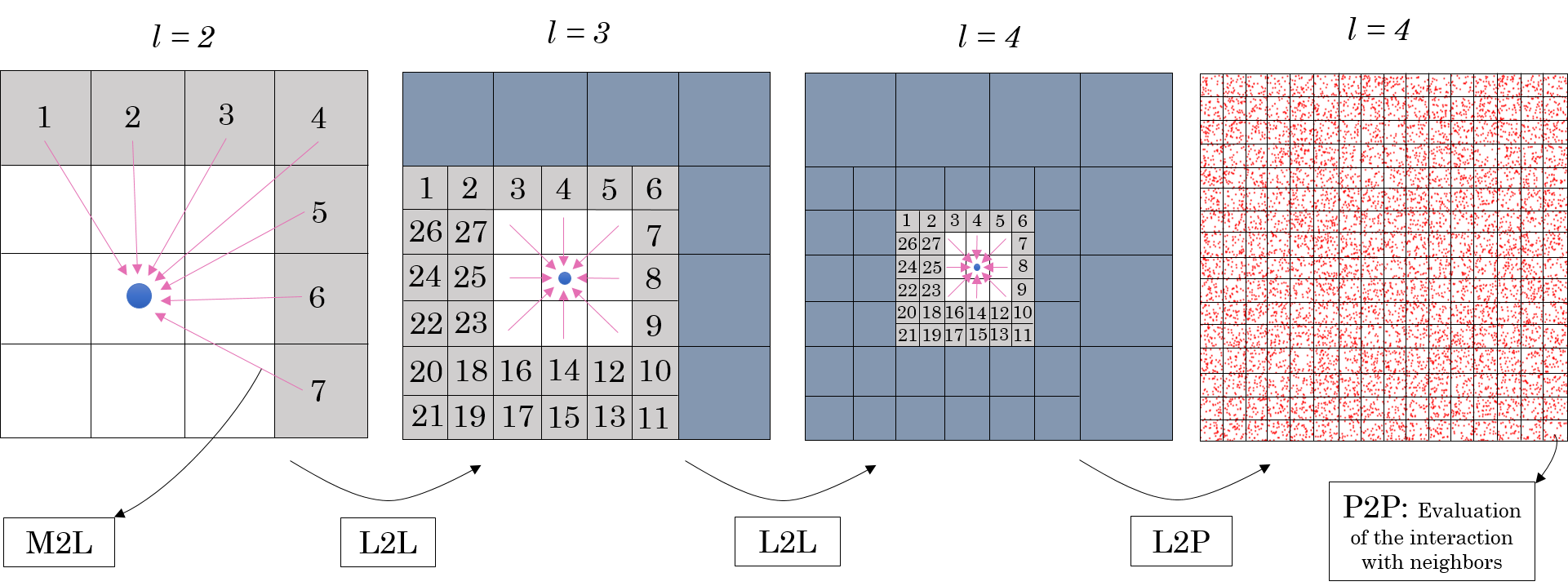}
		\caption{Downward pass} 
		\label{fig:DowneardPass}	
	\end{figure}
	Then, the \textit{downward pass} is performed for the final evaluation of the 
	potential at the influenced points.
	This aims at treating differently the contributions to the potential coming 
	from interaction list, far cells and neighbor cells.
	Figure \ref{fig:DowneardPass} shows the steps that are required.
	
	1)
	For a given cell $ I $, the contributions at the Chebyshev nodes $ 
	\bar{x}_{m}^I $ to the potential field $ f $ 
	coming from the interaction list of $ I $
	are calculated by the following \textit{multipole-to-local} (M2L) translation
	(let us denote with $ \mathrm{IntList}(I) $ the interaction list of $ I $) 
	\begin{align}
	( l = 0, 1, \ldots, k ) \quad g_m^I = \sum_{J \in \mathrm{IntList}(I)}
	\sum_{m'} W_{m'}^J G(\bar{x}_m^I,\bar{y}_{m'}^J), \quad m = 1,...,n  
	\label{m2l} \,.
	\end{align}
	
	2) In order to add 
	the contributions from far cells on the current cell $ I $, we travel the tree 
	downwards to get
	\begin{align}
	(l = 0)            \quad  f_m^I & = g_m^I \,, \quad m = 1,...,n \,, \\
	(l = 1,2,\ldots k) \quad  f_m^I & = g_m^I + \sum_{l'}f_{l'}^J 
	S_n(\bar{x}_m^I,\bar{x}_{l'}^J), \quad m = 1,...,n  \label{l2l} \,.
	\end{align}
	Equation \eqref{l2l} is called \textit{local-to-local} (L2L) interpolation 
	based on the parent cell $ J $ of $ I $.  
	
	3) Finally, the approximation $ \widetilde{f} $ of the potential is calculated 
	at each influenced point of $ I $ 
	by adding two terms:
	the local-to-point (L2P) term, for the interpolation of the influences $ f_m^I 
	$ coming from both interaction list and far cells,
	and the point-to-point (P2P) term, to account for the self and neighbor 
	interactions
	(we denote with $ \mathrm{Neighb}(I) $ the list of neighbors of $ I $)
	\begin{equation}
	( l = k ) \quad 
	\widetilde{f}(x_i) = 
	\sum_{m} f_m^I S_n(\bar{x}_m^I,x_i) + 
	\sum_{J \in \mathrm{Neighb}(I)} \sum_{y_j \in J} b(y_j) G(x_i,y_j)   \,.
	\end{equation}
	The procedure that was explained in this section 
	can be implemented within a BEM algorithm 
	such as the PDDM algorithm. This will be explained in the next section.
	
	%
	
	\section{PDDM with FMM (FMPDDM)} \label{sec:sol_proc}
	
	\begin{figure}[t]
		\centering
		\includegraphics[scale=0.3]{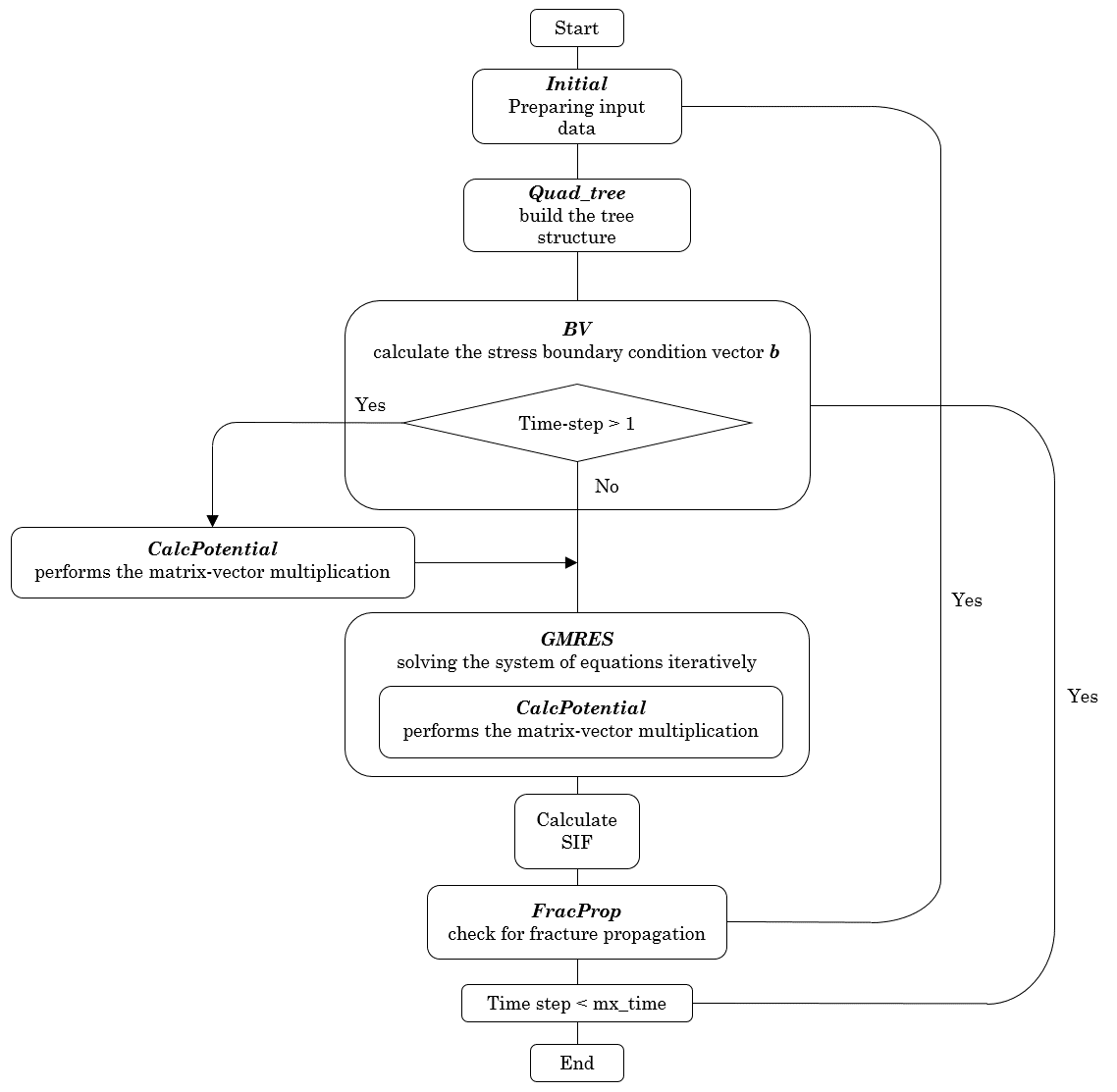}
		\caption{FMPDDM algorithm} 
		\label{fig:algorithm}	
	\end{figure}
	
	The FMPDDM algorithm that is used in this study is shown in Figure 
	\ref{fig:algorithm}. 
	As it can be seen, the first step to obtain $ D_s $, $D_n$ and $ D_q$ is to 
	provide 
	the spatial, temporal and material poroelastic properties, 
	and rock far-field stress and pore pressure,
	and the prescribed properties of the FMM model.
	Next, the quad-tree structure needs to be constructed. 
	It should be noted that the tree structure remains constant
	as long as the geometry is the same (i.e no fracture propagation occurs). 
	In the first time step, the right-hand side vector (i.e $\sigma_{s_0}$, 
	$\sigma_{n_0}$, and $P_{p_0}$) 
	is formed by transforming the far-field stresses to the face of each element 
	and by adding their effect to the hydraulic load on the elements. 
	For the time steps $ t = 2,...,\tau_h $ the solutions obtained from all of the 
	previous time steps 
	need to be subtracted from the known right-hand side vector 
	to account for the term with double summation in Equations 
	\eqref{eq:final_sigs}-\eqref{eq:final_pp}. 
	For this purpose, the matrix-vector multiplication between the coefficient 
	matrices 
	and the solutions of the previous time steps is performed by using FMM.
	Finally, a modified version of GMRES with FMM for matrix-vector multiplication 
	is used to iteratively solve for the problem unknowns.
	
	After obtaining the solutions, one may calculate the SIFs mode-I and II using 
	Equation \eqref{eq:SIF}. 
	Then, a fracture propagation criterion may be used to check to see whether the 
	propagation criterion is satisfied. 
	Consequently, an element will be added to the model if propagation happens, 
	and in that case the quad-tree needs to be constructed again. 
	In the case of no propagation, the solution of later time steps is obtained 
	by marching in time until reaching the final time step.

	\section{Numerical results} \label{sec:results}
	
	Some examples are presented here to demonstrate the application of the method. 
	The aim is to compare the accuracy and computational time of FMPDDM versus PDDM.
	The properties of Westerly granite used as input data for all the subsequent 
	examples
	are shown in Table \ref{table:rock_prop}.
	The properties that are required for the model are 
	shear modulus $ G $, 
	drained and undrained Poisson ratios $ \nu $ and $ \nu_d $, 
	Skemptson's coefficient $ B $,
	diffusivity coefficient $ c $,
	permeability $ k $, 
	and Biot's poroelasticity coefficient, $ \alpha $. 
	
	\begin{table}[t]
		\centering
		\begin{tabular}{ | l | c | c |c | c |c | c |c |}
			\hline
			\textbf{Rock type} & $ G $ & $\nu$ & $\nu_u$ & $B$ & $ c $ & $k$ & $ 
			\alpha $ \\  
			& $ GPa $ & $-$ & $ - $ & $-$ & $ m^2/s $ & $ m^2$ & $ - $ \\ 		
			\hline
			\hline
			Westerly granite  & $15$ & $0.25$ & $0.331 $ & 
			$0.81$ & $6.15 \cdot 10^{-5}$ & $4 \cdot 
			10^{-19}$  & 
			$0.449$  \\ \hline
		\end{tabular}
		\caption{Rock poroelastic properties used in this study (adopted from 
		\cite{cheng2016poroelasticity}).}
		\label{table:rock_prop}
	\end{table}

	\subsection{Dependence on the number of fracture elements}
	
	\begin{figure}[t]
		\centering
		\includegraphics[scale=0.2]{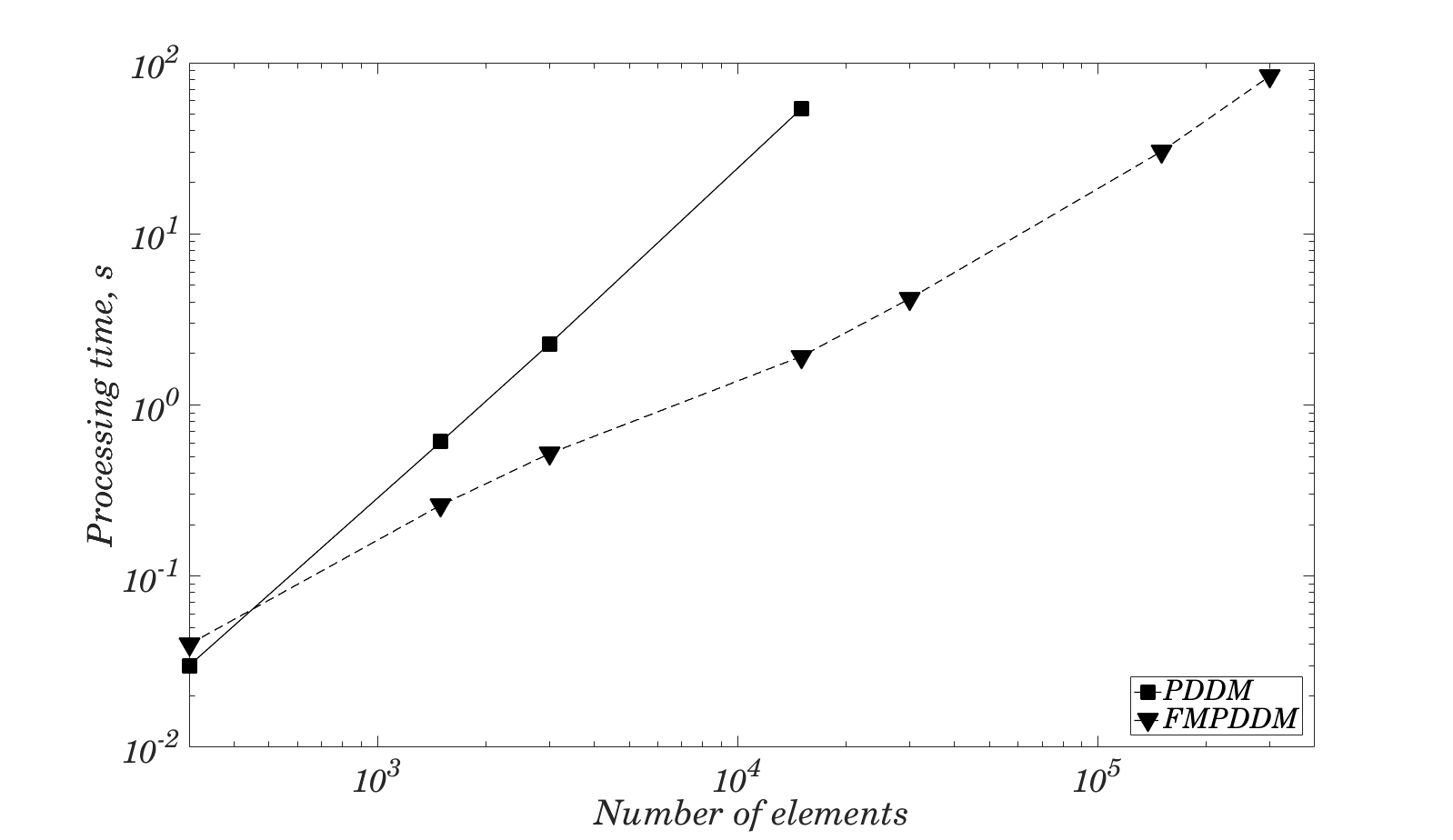}
		\caption{CPU processing time vs. number of elements for one time step} 
		\label{fig:CPU_one_step}	
	\end{figure}

	First, the processing times of the two methods are compared as a function of 
	the number of boundary elements. 
	The comparison for one time step of the solution is shown in Figure 
	\ref{fig:CPU_one_step}.
	Initially PDDM performs better in terms of processing time,
	but, as the number of elements increases, 
	it takes more time for the PDDM method to solve a problem with the same number 
	of elements in one time step.
	By further increasing the number of elements, the slope of PDDM becomes equal 
	to two,
	while the FMPDDM has the slope of one.
	This indicates that PDDM has the complexity of $ O (N^2) $, 
	while FMPDDM has the complexity of $ O(N) $.
	
	\begin{figure}[t]
		\centering
		\includegraphics[scale=0.2]{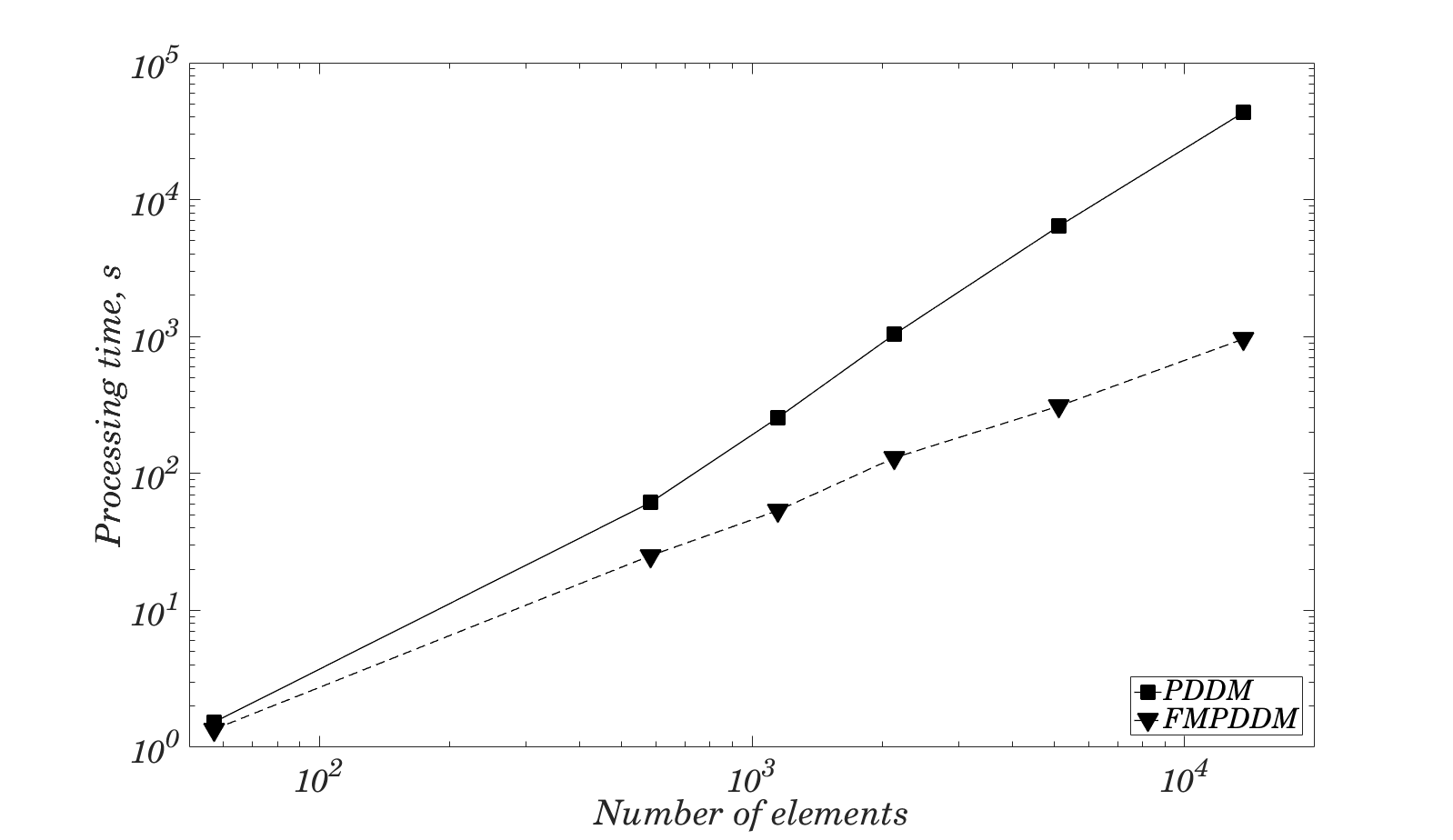}
		\caption{CPU processing time vs. number of elements for 10 time steps} 
		\label{fig:CPU_10_step}	
	\end{figure}

	Figure \ref{fig:CPU_10_step} shows the processing times for 10 time steps with 
	the same time discretization
	as a function of the number of elements. 
	As required by the time marching method,
	for all time steps other than the first time step,
	it is necessary to subtract all the previous solutions from the current time 
	step boundary condition 
	(Equations \eqref{eq:final_sigs}-\eqref{eq:final_pp}). 
	This creates extra steps with matrix-vector multiplications, 
	increasing the processing time even further. 
	As it may be seen in the figure, as the number of elements increases,
	FMPDDM performs much better than PDDM in terms of processing time. 
	For example, for $ 13682 $ elements it was observed that the processing time of 
	the PDDM method was about $ 43340 $ seconds, 
	while for the same problem it took only $ 960 $ seconds 
	for the FMPDDM to calculate the solution. 
	This is a huge difference in processing time ($\sim 40$ times less). 
	Also, a better performance is expected for greater number of elements.
	
	\subsection{Accuracy of the solution}

	Another important aspect to investigate for FMDDM is to see how the accuracy of 
	the solution is preserved.
	In order to investigate this feature, two examples are investigated.
	In the first example, a hydraulic fracturing problem in a horizontal well is 
	studied. 
	In the second example, the problem of randomly distributed pressurized 
	fractures is presented. 
	For each case, shear and normal displacements of a certain element 
	as well as its stress intensity factors  $ K_{I} $ and $ K_{II} $ are 
	calculated. 
	For each case, the accuracy is discussed along with the factors that may be 
	used to improve it.

	\subsubsection{Transverse Hydraulic Fractures in a Cluster}
	
\begin{figure}[H]
		\centering
		\subfloat[]{
			\includegraphics[width=0.4\textwidth]{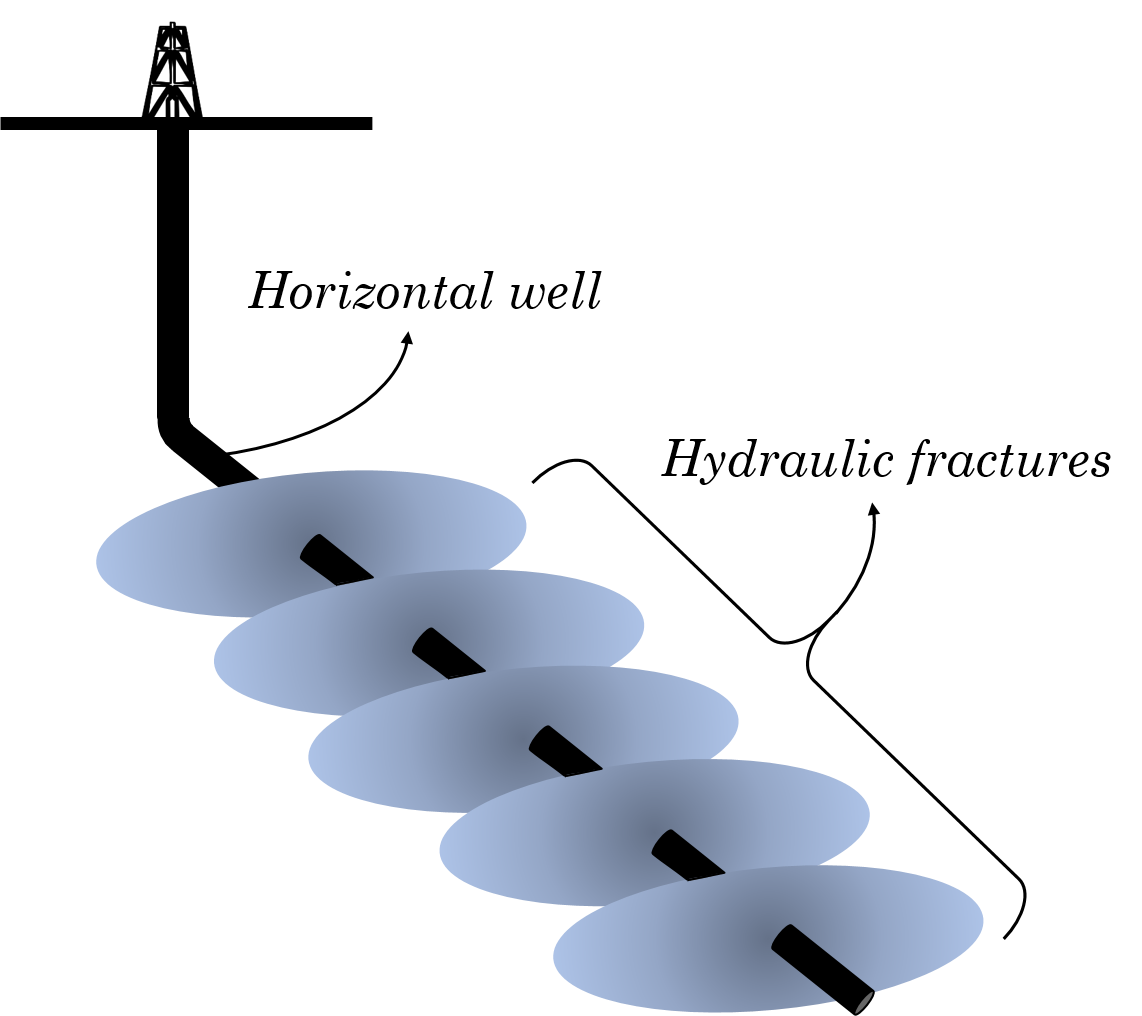}
			\label{fig:hf_3d_view}
		}
		\subfloat[]{
			\includegraphics[width=0.4\textwidth]{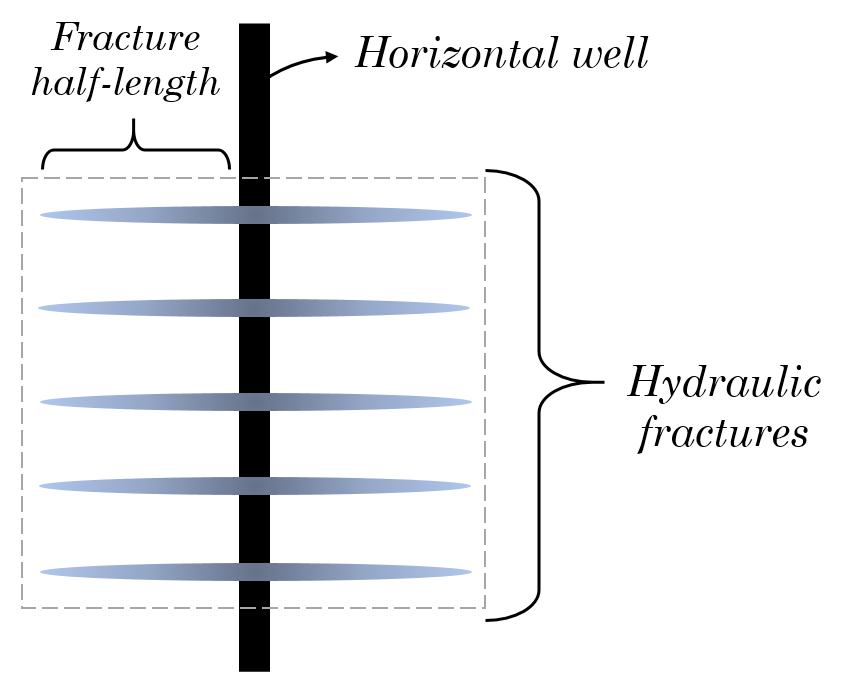}
			\label{fig:hf_top_view}
		}	
		\caption{Hydraulic fracturing in a horizontal well; a) 3D view, 
			b) top view}
		\label{fig:multiple_fractures_horizontal_well}
	\end{figure}

	It is a common practice in the oil industry to drill a well horizontally
	and initiate multiple hydraulic fractures to increase the hydrocarbon 
	production.
	In most of the ultra-tight reservoirs, hydrocarbon recovery is impossible 
	without hydraulic fractures.
	Figure \ref{fig:hf_3d_view} shows a schematic of the hydraulic fracture process 
	in a horizontal wellbore. 
	Usually, the well is drilled parallel to the minimum horizontal stress direction
	and multiple sections of the wellbore are isolated and perforated from the 
	``toe'' of the horizontal well to its ``heel''.
	This causes hydraulic fractures to initiate and propagate orthogonal to the 
	wellbore 
	(i.e. parallel to the direction of maximum horizontal compressional stress).
	Consequently, new interface areas are created inside the reservoir 
	that cause an increase in the production of hydrocarbon from wellbore.   
	
	In order to show the accuracy of the proposed algorithm,
	the solution of five parallel hydraulic fractures in a horizontal well is 
	studied. 
	Figure \ref{fig:hf_top_view} shows the geometry of the problem that is used for 
	this purpose. 
	Five parallel hydraulic fractures, each having $ 400 \, m $ half-length 
	(hydraulic length) are assumed.
	Each fracture is discretized using $ 200 $, $ 300 $, and $ 400 $ constant DDM 
	elements.

	The boundary condition and the other inputs of the model are summarized in 
	Table \ref{table:hydraulic_frac}. 
	For this problem, we assume that the directions of maximum and minimum 
	horizontal stresses are $x$ and $y$, respectively. 
	The spacing of each fracture is also chosen to be $ 30 \; m $.
	The fluid is injected into the wellbore with $ 58\; MPa $ and the problem is 
	solved for the first time step equal to $ 100\;s $.

	\begin{table}[H]
		\centering
		\begin{tabular}{ | l | c |}
			\hline
			Maximum horizontal stress, $\sigma_H$ & $58.60 \text{ MPa}$ \\ \hline
			Minimum horizontal stress, $\sigma_h$ & $55.15 \text{ MPa}$ \\ \hline
			Initial reservoir pore pressure, $p_0$ & $48.26 \text{ MPa}$ \\ \hline 
			Injection pressure, $p_i $ & $58 \text{ MPa}$ \\ \hline 
			Fracture spacing	&	 $30 \text{m}$ \\ \hline
			Fracture half-length, $a$	& $400 \text{m}$ \\ \hline
		\end{tabular}
		\caption{Model configuration for hydraulic fractures in a horizontal 
		wellbore}
		\label{table:hydraulic_frac}
	\end{table}
	
	In order to compare the accuracy of the model, we have chosen four elements 
	and compared the shear and normal displacements on these elements using FMPDDM 
	and PDDM. 
	The center elements and tip elements of the inner and outer fracture
	(third and fifth fracture from the top in Figure 
	\ref{fig:hf_top_view}).
	are the elements that we chose for this purpose.
	The number of Chebyshev nodes is set to 6 for this example. 
	Table \ref{tab:Dn_inner} shows the normal displacement of the center element 
	and tip element 
	of the inner fracture for different discretizations of $ 200 $, $ 300 $, and $ 
	500 $
	up to four decimal places for each fracture respectively.

	\begin{table}[H]
		\centering
		\begin{tabular}{ | l | c | c |c | c |c | c |}
			\hline
			$D_n$ & \multicolumn{3}{c|}{\textit{Fracture center}} & 
			\multicolumn{3}{c|}{\textit{Fracture tip}} \\ \hline 
			NE	&	$1000$	&	$1500$	&	$2000$	&	$1000$	&	$1500$ & $2000$ 
			\\  
			\hline
			Max. level & 2 & 3 & 3 & 2 & 3 & 3   \\  
			\hline
			PDDM	& $0.0199$	& $0.0184$ & $0.0170$ & $0.0049$ &	$0.0040$	&  
			$0.0034$ \\ \hline
			FMPDDM	& $0.0193$	& $0.0179$ & $0.0166$ & $0.0047$ &	$0.0039$	&  
			$0.0033$ \\ \hline
			error (\%)	& $-3$	& $-2.7$ & $-2.3$ & $-4$ &	$2.5$	&  $-2.9$ \\ 
			\hline
		\end{tabular}
		\caption{Normal displacement of the inner fracture surface for different 
		number of elements (NE) after the first time step}
		\label{tab:Dn_inner}
	\end{table}

	As it may be seen in the Table \ref{tab:Dn_inner}, the normal opening of the 
	fracture at the center is higher as expected. 
	Also, the maximum error that is reported for normal displacements is smaller 
	than 4\% for all cases. 
	Also, it should be noted that the results that are obtained by FMPDDM are 
	smaller than the results that are obtained by PDDM. 
	Shear displacements at the center and tip of the fracture were smaller than the 
	tolerance of GMRES ($ = 10^{-6} $).
	This was expected because having any shear displacement on the elements in the 
	middle fracture 
	causes a mixed mode (mode I + II) fracture propagation. 
	This is not the case for such an arrangement as was shown previously. 
	Therefore, we didn't report the shear displacement of the center fracture for 
	this case. 
	Similar calculation may be done for the outer fracture. 
	Table \ref{tab:Dn_outer} shows the normal displacement of the outer fracture.

	\begin{table}[H]
		\centering
		\begin{tabular}{ | l | c | c |c | c |c | c |}
			\hline
			$D_n$ & \multicolumn{3}{c|}{\textit{Fracture center}} & 
			\multicolumn{3}{c|}{\textit{Fracture tip}} \\ \hline 
			NE	&	$1000$	&	$1500$	&	$2000$	&	$1000$	&	$1500$ & $2000$ 
			\\  
			\hline
			PDDM	& $0.0198$	& $0.0183$ & $0.0170$ & $0.0057$ &	$0.0044$	&  
			$0.0037$ \\ \hline
			FMPDDM	& $0.0192$	& $0.0178$ & $0.0166$ & $0.0056$ &	$0.0043$	&  
			$0.0036$ \\ \hline
			error (\%)	& $-3$	& $-2.7$ & $-2.3$ & $-1.75$ &	$-2.27$	&  $-2.7$ 
			\\ \hline
		\end{tabular}
		\caption{Normal displacement of the outer fracture surface for different 
		number of elements (NE) after the first time step}
		\label{tab:Dn_outer}
	\end{table}
	
	Similar to what was observed for the inner fracture normal displacement, 
	the error of the normal displacement on the outer fracture was smaller than 
	4\%. 
	Also, as expected, both shear and normal openings of the outer fracture were 
	bigger compared to the inner fracture. 
	This is because the stress shadow that is created by other outer fractures 
	impedes the inner fracture from opening.
	Next, shear displacement of the outer fracture is presenter in Table 
	\ref{tab:Ds_outer}. 
	For the same reason that explained in the case of inner fracture shear 
	displacement, 
	only the shear displacement of one of the outer fracture tips is presented 
	here.

	\begin{table}[H]
		\centering
		\begin{tabular}{ | l | c | c |c | c |c | c |}
			\hline
			$D_s$  & 
			\multicolumn{3}{c|}{\textit{Fracture tip}} \\ \hline 
			NE	&	$1000$	&	$1500$ & $2000$ \\  
			\hline
			PDDM	&  $0.0024$ &	$0.0019$	& $0.0015$  \\ \hline
			FMPDDM	&  $0.0023$ &	$0.0018$	& $0.0015$  \\ \hline
			error (\%)	&  $-4.16$ &	$-5$	&  $0$ \\ \hline
		\end{tabular}
		\caption{Shear displacement of the outer fracture surface for different 
		number of elements (NE) after the first time step 
			(Note: the shear displacement of the fracture center was less than 
			machine precision).}
		\label{tab:Ds_outer}
	\end{table}
	
	As shown in Table \ref{tab:Ds_outer}, the same error was also observed for the 
	shear displacement of the outer fracture. 
	Also, it was expected that because of stress shadowing effect,
	some shear will be observed at the tip of the outer fracture since the outer 
	fractures tend to reorient away from the middle fractures.
	As shown in this example using FMPDDM will give acceptable results with smaller 
	computation time. 
	The error was smaller than 4\% in all of the cases.
	As an example, the processing time required for FMPDDM to solve the $ 2000 $  
	element case ($ 400 $ element per fracture) 
	was one third of the time required for PDDM to solve the same problem. 
	Next, we present a case of randomly distributed fractures and do the same 
	exercise 
	while giving different angles to the distributed fractures. 
	
	\subsubsection{Randomly Distributed Pressurized Fractures}
	
	Underground rocks are discontinuous media filled with natural fractures.
	Thus, it is necessary to include them in any hydraulic fracturing study.
	In this section, randomly distributed pressurized fractures are studied
	to analyze the accuracy of the developed FMPDDM.
	Here, we assume that all fractures are pressurized. 
	Similar to the previous section, we choose one fracture here and calculate the 
	normal and shear displacements of both tips.
	Then, we calculate mode I and II stress intensity factors that are crucial for 
	any fracture propagation study.  
	
\begin{figure}[!h]
		\centering
		\subfloat[$15^o$]{
			\includegraphics[width=0.4\textwidth]{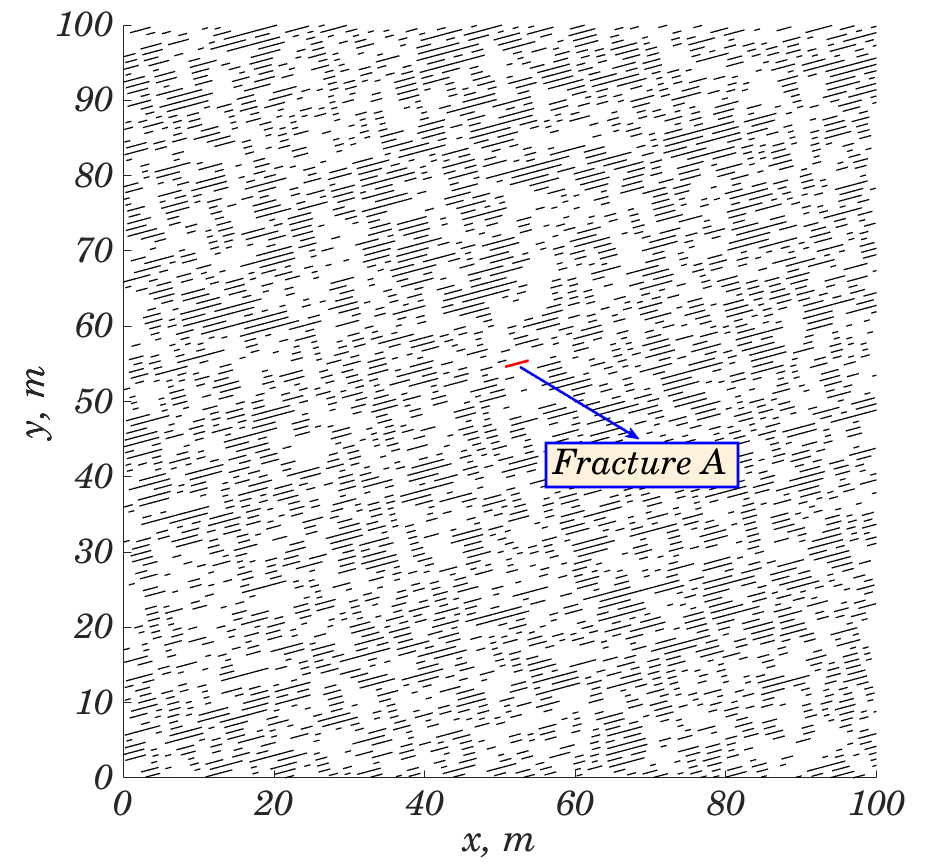}
			\label{fig:random_15}
		}
		~
		\subfloat[$45^o$]{
			\includegraphics[width=0.4\textwidth]{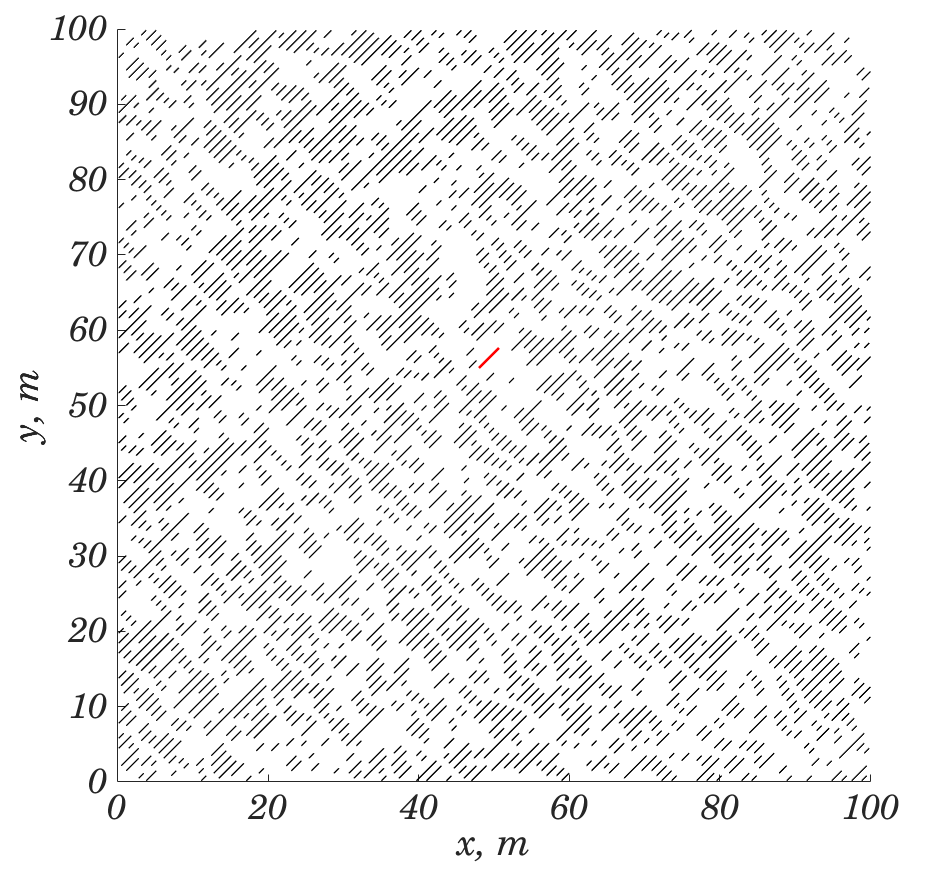}
			\label{fig:random_35}
		}
	    \quad
		\subfloat[$75^o$]{
			\includegraphics[width=0.4\textwidth]{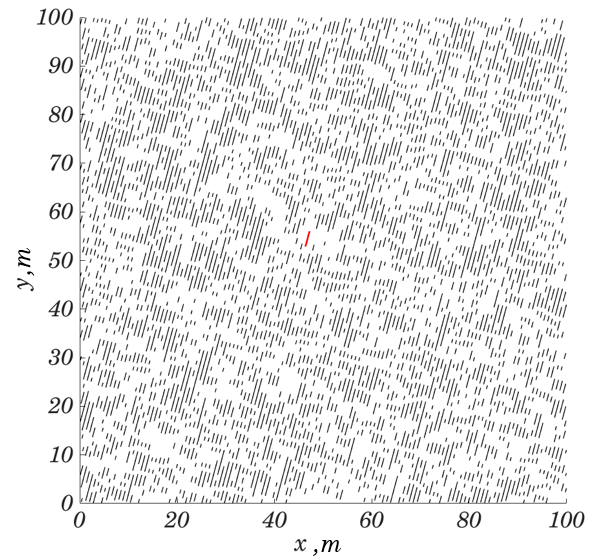}
			\label{fig:random_55}
		}	
		\caption{Randomly distributed fractures (angles are with respect to maximum 
		($x$ axis) horizontal stress).}
		\label{fig:randomly_distributed_fractures}
	\end{figure}

	Figure \ref{fig:randomly_distributed_fractures} shows the geometry of the 
	problem that is studied in this section.
	The fracture that is selected for the displacement calculation is shown on the 
	figure (fracture \textit{\textit{A}}). 
	Fracture \textit{A} has five elements. 
	It is also assumed that all fractures are parallel to each other and three 
	cases 
	with different angles of $ 15^o $, $ 45^o $, $ 75^o $ with respect to the $ x $ 
	axis are considered. 
	The direction of maximum and minimum horizontal stresses in all cases are $ x $ 
	and $ y $, respectively.

	Table \ref{table:random_frac} represents the boundary conditions and the other 
	input variables of the model. 
	The injection pressure in this case is $ 60\;MPa $ to make sure all the 
	fractures are open. 
	Time is discretized into one time step equal to $ 10 $ seconds.
	Also, to increase the accuracy of FMPDDM, Chebyshev polynomials of degrees 
	three (NCheb3) and six (NCheb6) are chosen.
	The processing times of  PDDM, FMPDDM-NCheb6, and FMPDDM-NCheb3 observed as 
	8050, 421, and 120 seconds, respectively. 
	Therefore, the processing time can be cut by almost 67 times.

	\begin{table}[H]
		\centering
		\begin{tabular}{ | l | c |}
			\hline
			Maximum horizontal stress, $\sigma_H$ & $58.60 \text{ MPa}$ \\ \hline
			Minimum horizontal stress, $\sigma_h$ & $55.15 \text{ MPa}$ \\ \hline
			Reservoir pore pressure, $p_r$ & $48.26 \text{ MPa}$ \\ \hline 
			Pressure inside the fracture, $p_f $ & $60 \text{ MPa}$ \\ \hline 
			Number of fractures	& $4335$  \\ \hline
			Number of elements	& $20332$ \\ \hline
			Number of fracture \textit{A} elements	& $5$ \\ \hline
		\end{tabular}
		\caption{Input parameters of the randomly distributed pressurized 
		fractures.}
		\label{table:random_frac}
	\end{table}


	Table \ref{tab:DS_FracA} shows the shear displacements at the right and left 
	tips of fracture $ A $ for different angles.
	Comparing PDDM, FMPDDM-NCheb3, and FMPDDM-NCheb6 shows that, 
	although using the FMPDDM-NCheb3 requires only two minutes compared with more 
	than two hours required with PDDM,
	it causes a small error in the solution. 
	The error may be reduced by using FMPDDM-NCheb6, 
	but it increases the computational time to $ 7 $ minutes, which is still much 
	less than the run time required for PDDM.

	\begin{table}[H]
		\centering
		\begin{tabular}{ | l | c |c | c | c |c | c|}
			\hline
			\textbf{$D_s$}  & 
			\multicolumn{3}{c|}{\textit{Fracture left tip $\times 10 ^{-5}$ m}} & 
			\multicolumn{3}{c|}{\textit{Fracture right tip $\times 10 ^{-5}$ m}}\\ 
			\hline 
			Angle & $15^o $ &  $45^o$  & $75^o$ & $15^o $ & $45^o $  & $75^o$ \\    
			\hline
			Conventional PDDM              & $2.0592$ & $7.2157$ & $2.6829$ & 
			$4.6915$ & $8.1299$ & $2.6438$ \\ \hline
			FMPDDM-NCheb3    & $2.0545$ & $7.2050$ & $2.6795$ & $4.6736$ & $8.1307$ 
			& $2.6478$ \\ \hline
			FMPDDM-NCheb6    & $2.0594$ & $7.2155$ & $2.6829$ & $4.6915$ & $8.1298$ 
			& $2.6438$ \\ \hline
		\end{tabular}
		\caption{Shear displacement of the tips of fracture $A$ }
		\label{tab:DS_FracA}
	\end{table}
	
	The same calculation for the normal displacement of fracture $ A $ tips is 
	presented in Table \ref{tab:Dn_FracA}. 
	Similar to shear displacement, the same error was observed for this case. 
	In both cases, the relative error in calculating the displacement is in the 
	order of $10^{-3}$, 
	and FMPDDM-NCheb3 error is four times bigger than FMPDDM-NCheb6. 
	The error for both cases is negligible, and as it will be discussed in the next 
	section,
	it will not change the SIF calculation, which is the key part for calculating 
	the fracture propagation.
	
	\begin{table}[H]
		\centering
		\begin{tabular}{ | l | c |c | c | c |c | c|}
			\hline
			\textbf{$D_n$}  & 
			\multicolumn{3}{c|}{\textit{Fracture left tip $\times 10 ^{-5}$ m}} & 
			\multicolumn{3}{c|}{\textit{Fracture right tip $\times 10 ^{-5}$ m}}\\ 
			\hline 
			Angle & $15^o $ &  $45^o$  & $75^o$ & $15^o $ & $45^o $  & $75^o$ \\    
			\hline
			Conventional PDDM             & $2.2694$ & $11.733$ & $1.2616$ & 
			$12.227$ & $11.443$ & $0.8101$ \\ \hline
			FMPDDM-NCheb3    & $2.3158$ & $11.704$ & $1.2612$ & $12.246$ & $11.460$ 
			& $0.8113$ \\ \hline
			FMPDDM-NCheb6    & $2.2697$ & $11.733$ & $1.2616$ & $12.227$ & $11.443$ 
			& $0.8101$ \\ \hline
		\end{tabular}
		\caption{Normal displacement of fracture $A$}
		\label{tab:Dn_FracA}
	\end{table}
	
	The accuracy and speed of the PDDM and two FMPDDMs via two examples were 
	presented in this section. 
	The relative errors for the case of parallel fractures was calculated to be 
	less than $ 4 $\%. 
	The reason for the big differences in the error is that in the first example,
	elements were not evenly distributed on a plane, 
	whereas in the second example points where distributed on the plane. 
	In the next section, the stress intensity factor for the right tip of fracture 
	$ A $ is calculated 
	using FMPDDM with two different degrees, and the results are compared to the 
	stress intensity factor 
	calculated using the PDDM.

	
	The stress intensity factor calculation 
	is the key part of any fracture propagation investigation. 
	Any error in the calculation of SIF results in erroneous determination of both 
	angle 
	(if considering mixed mode I+II) and the moment of fracture extension. 
	In order to investigate the accuracy of the developed FMPDDM model, 
	the SIF calculation on the right tip of fracture $ A $ is conducted. 
	Table \ref{tab:SIF_FracA} shows the result of this calculation. 
	The error of FMPDDM-NCheb6 is less than $  10^{-4} $ for this problem.   
	
	\begin{table}[H]
		\centering
		\begin{tabular}{ | l | c | c |c | c |c | c |}
			\hline
			{} & \multicolumn{3}{c|}{\textit{$SIF_I$, $ MPa \sqrt{m}$ }} & 
			\multicolumn{3}{c|}{\textit{$SIF_{II}$, $ MPa\sqrt{m}$  }} \\ \hline 
			Angle & $ 15^o $ & $45^o$ & $75^o$ & $15^o $ & $ 45^o $  & $75^o$ \\    
			\hline
			Conventional PDDM    & $3.0304$ & $2.8362$ & $0.20078$ & $1.1628$ & 
			$2.0150$ & $0.65527$  \\ \hline
			FMPDDM-NCheb3        & $3.0352$ & $2.8405$ & $0.20109$ & $1.1584$ & 
			$2.0152$ & $0.65627$  \\ \hline
			FMPDDM-NCheb6        & $3.0305$ & $2.8362$ & $0.20078$ & $1.1628$ & 
			$2.0150$ & $0.65527$  \\ \hline
		\end{tabular}
		\caption{Stress intensity factor mode II ($K_{II}$) at the right tip of 
		fracture $A$}
		\label{tab:SIF_FracA}
	\end{table}
	
	The errors of SIF calculated by FMPDDM using both polynomials of degree 3 and 6 
	are small enough  
	and do not cause an error in the calculation of the angle and the moment of 
	fracture extension. 
	Fracture propagation is subject of future investigations.    
	
	\section{Conclusions}
	
	In this study, a kernel-independent fully-poroelastic fast multipole 
	displacement discontinuity model was developed to study hydraulic fracturing 
	problems. The FMM uses Chebyshev polynomials for the kernel expansion.
	Dependence of the model performance to number of elements is compared with the 
	conventional fully-poroelastic DDM for one kernel matrix-vector multiplication, 
	and a complete hydraulic fracture problem including 10 time steps. Also, in 
	order to compare the accuracy of the developed method, 
	its performance was compared with a conventional displacement discontinuity 
	method using two examples.
	In the first example the injection of the fracturing fluid into five parallel 
	hydraulic fractures in a horizontal well was studied. 
	It was observed that the maximum error in the shear and normal displacements of 
	the fracture surfaces at the wellbore and fracture tips 
	is less than $ 5 $\%.
	This error may be further minimized by choosing higher order Chebyshev 
	polynomials. 
	Generally, choosing higher order polynomials requires more computational time.
	In the second example, the problem of randomly pressurized fractures was 
	considered. 
	The error in this case was about $ 0.5 $\%. 
	The difference between the error in the two examples is due to the distribution 
	of the elements 
	in the plane in the second example.
	In the first example, points were on the fractures,
	while in the second example points were distributed randomly in the plane.  
	The results show that hydraulic fracture problems may be solved up to 
	approximately 70 times faster 
	by incorporating the fast multipole method 
	into the conventional poroelastic displacement discontinuity method 
	for a case of 2000 elements and one time step. 
	Increasing the number of elements or time steps will make the ratio even 
	bigger. 
	We also showed that the error introduced by FMM in the solution is small enough 
	that it has no significant effect on the SIF calculation, 
	while decreasing the processing time. 
	Therefore, fracture propagation angle and its occurrence moment can be 
	estimated 
	with high degree of confidence with FMPDDM.
	Although using higher order polynomials will affect the processing time and 
	accuracy of the solution with FMPDDM, 
	there are some other key factors such as iterations number in GMRES, degree of 
	Chebyshev polynomials, tolerance of the GMRES, 
	and distribution of the elements in a plane that also affect the results. 
	In the next part of this study, the propagation of fractures will be studied 
	using the SIFs 
	that were calculated in this paper.

	\bibliographystyle{spbasic}  
	\bibliography{biblio} 
		
\end{document}